\RequirePackage{fix-cm}
\documentclass{article}
\usepackage{graphicx}
\usepackage{amsmath,amssymb,,amsthm}
\usepackage{graphicx}
\usepackage{mathrsfs}

\newtheorem{theorem}{Theorem}
\newtheorem{lemma}{Lemma}
\newtheorem{corollary}[lemma]{Corollary}
\theoremstyle{remark}
\newtheorem{remark}{Remark}
\newtheorem{definition}[lemma]{Definition}

\def \IR {\mathbb{R}}
\def \Th {{\cal T}_h}
\def\proof{\noindent{\sl Proof: }}
\def\qed{$\hfill\square$}
\def\HO0{{H^1_0(\Omega)}}
\def\HO{{H^1(\Omega)}}
\def\LO{{L^2(\Omega)}}
\def\LG{{L^2(\Gamma)}}
\def\LoL{{L^2(0,L)}}
\def\HoL{{H^1(0,L)}}
\def\HK{{H^1(K)}}
\def\HdK{{H^2(K)}}
\def\Wi{{W^{2,\infty}(\Omega)}}
\def\grad{\nabla}
\def\vhK{{v_h^{K}}}
\def\pd{\partial}
\def\di{{2,\infty}}
\def \d {{\rm d}}
\def\dx{\,\d x }
\def\dS{\,\d S }
\def \n{n}
\def \bmu{\mathbf{\mu}}
\def \B{\mathcal{B}}
\def\diver{{\rm div}}
\def\diam{{\rm diam\,}}

\def\UG{U_\Gamma} 
\def\uG{u_\Gamma} 
\def\wp{w^\prime}
\def\Pih{\Pi_h}
\def\tu{\tilde{u}}
\def\T{{\mathrm{T}}}
\def\vphi{\varphi}
\def\rK{{r_K}}
\def\supp{\mathrm{supp}\, }
\def\talpha{{\tilde{\alpha}}}

\def\mcA{\mathcal{A}}

\def\mcC{\mathcal{C}}
\def\BB{{\mathbb{B}}\, }
\def\mC{{\mathscr{C}}}
\def\vhC{v_h^{\mathscr{C}}}
\def\lmC{{\!_\mC}}
\def\tK{{\tilde{K}}}

\begin{document}

\title{On necessary and sufficient conditions for finite element convergence\thanks{This work is a part of the research project P201/13/00522S of the Czech Science Foundation. V. Ku\v cera is currently a Fulbright visiting scholar at Brown University, Providence, RI, USA, supported by the J. William Fulbright Commission in the Czech Republic.}}


\author{V\' aclav Ku\v cera \\ \\Faculty of Mathematics and Physics, Charles University in Prague}



\maketitle

\begin{abstract}
In this paper we derive a necessary condition for finite element method (FEM) convergence in $H^1(\Omega)$ as well as generalize known sufficient conditions. We deal with the piecewise linear conforming FEM on triangular meshes for Poisson's problem in 2D. In the first part, we prove a necessary condition on the mesh geometry for $O(h^\alpha)$ convergence in the $H^1(\Omega)$-seminorm with $\alpha\in[0,1]$. We prove that certain structures, bands consisting of neighboring degenerating elements forming an alternating pattern cannot be too long with respect to $h$. This is a generalization of the Babu\v ska-Aziz counterexample and represents the first nontrivial necessary condition for finite element convergence. Using this condition we construct several counterexamples to various convergence rates in the FEM. In the second part, we generalize the maximum angle and circumradius conditions for $O(h^\alpha)$ convergence. We prove that the triangulations can contain many elements violating these conditions as long as their maximum angle vertexes are sufficiently far from other degenerating elements or they form clusters of sufficiently small size. While a necessary and sufficient condition for $O(h^\alpha)$ convergence in $H^1(\Omega)$ remains unknown, the gap between the derived conditions is small in special cases.

\end{abstract}

\section{Introduction}
\label{sec:intro}
The finite element method (FEM) is perhaps the most popular and important general  numerical method for the solution of partial differential equations. In its classical, simplest form, the space of piecewise linear, globally continuous functions on a given partition (triangulation) $\Th$ of the spatial domain $\Omega$ is used along with a weak formulation of the equation. In our case, we will be concerned with Poisson's problem in 2D.

Much work has been devoted to the a priori error analysis of the FEM. Namely, the question arises, what is the necessary and sufficient condition for the convergence of the method when the meshes $\Th$ are refined, i.e. we have a system of triangulation $\{\Th\}_{h\in(0,h_0)}$. In the simplest case, we are interested in the energy norm, i.e. $H^1(\Omega)$-estimates of the error
\begin{equation}
|u-U|_1\leq C(u)h,
\label{sec:intro:est}
\end{equation}
where $u$ and $U$ are the exact and approximate solutions, respectively, $h$ is the maximal diameter of elements from $\Th$ and $C(u)$ is a constant independent of $h$. Typically, one is interested in deriving (\ref{sec:intro:est}) for some larger class of functions, e.g. for all $u\in H^2(\Omega)$.

Historically, the first sufficient condition for (\ref{sec:intro:est}) to hold is the so-called \emph{minimum angle condition} derived independently in \cite{Zenisek}, \cite{Zlamal}. This condition states that there should exist a constant $\gamma_0$ such that for any triangulation $\Th, h\in(0,h_0),$ and any triangle $K\in \Th$ we have 
\begin{equation}
0<\gamma_0\leq\gamma_K,
\label{sec:intro:min}
\end{equation}
where $\gamma_K$ is the minimum angle of $K$. This condition was later weakened independently by \cite{Babuska-Aziz}, \cite{Barnhill} and \cite{Jamet} to the \emph{maximum angle condition}:
there exists a constant $\alpha_0$ such that for any triangulation $\Th, h\in(0,h_0),$ and any triangle $K\in \Th$ we have 
\begin{equation}
\alpha_K\leq\alpha_0<\pi,
\label{sec:intro:max}
\end{equation}
where $\alpha_K$ is the maximum angle of $K$. Finally, a condition for convergence was recently derived in \cite{Kobayashi} and \cite{Rand}, the \emph{circumradius condition}: Let 
\begin{equation}
\max_{K\in\Th}R_K\to 0,
\label{sec:intro:circ}
\end{equation}
where $R_K$ is the circumradius of $K$. Then the FEM converges in $H^1(\Omega)$. We note that all these conditions are derived by taking the piecewise linear Lagrange interpolation $\Pi_h u$ of $u$ in C\' ea's lemma.

Conditions (\ref{sec:intro:min}), (\ref{sec:intro:max}) are sufficient conditions for $O(h)$ convergence. It was shown in \cite{Hannukainen} that the maximum angle condition is not necessary for (\ref{sec:intro:est}) to hold. The argument is simple and can be directly extended to the circumradius condition. We take a system of triangulations satisfying (\ref{sec:intro:max}) -- thus exhibiting $O(h)$ convergence -- and refine each $K\in\Th$ arbitrarily to obtain $\tilde\Th$. Thus $\tilde\Th$ can contain an arbitrary amount of arbitrarily bad `degenerating' elements, however since $\tilde\Th$ is a refinement of $\Th$, it also exhibits $O(h)$ convergence.

Since (\ref{sec:intro:min})--(\ref{sec:intro:max}) are only sufficient and not necessary, the question is what is a necessary and sufficient condition for (\ref{sec:intro:est}) to hold. The only step in this direction is the Babu\v ska-Aziz counterexample of \cite{Babuska-Aziz} consisting of a regular triangulation of the square consisting \emph{only} of `degenerating' elements violating the maximum angle condition (except for several cut-off elements adjoining to the boundary $\partial\Omega$), cf. Figure \ref{fig:Babuska-Aziz}. Recently a more detailed, optimal analysis of this counterexample was performed in \cite{Oswald}. By controlling the speed of degeneration of the aspect ration of the elements, one can produce a counterexample to $O(h)$ convergence and even an example of nonconvergence of the FEM. From this one counterexample, Babu\v ska and Aziz conclude that the maximum angle condition is \emph{essential} for $O(h)$ convergence, carefully avoiding the word necessary. 

In this paper, we will be concerned with a more general version of (\ref{sec:intro:est}), namely $O(h^\alpha)$ estimates of the form
\begin{equation}
|u-U|_1\leq C(u)h^\alpha,
\label{sec:intro:est_alpha}
\end{equation}
for $\alpha\in[0,1]$. We will derive a necessary condition on the mesh geometry and also a new sufficient condition for (\ref{sec:intro:est_alpha}) to hold. To the author's best knowledge, the derived necessary condition is the first nontrivial necessary condition for $O(h)$ or any other convergence of the finite element method, apart from the trivial necessary condition $h\to 0$. We note that the results of this paper are in fact about the approximation properties of the piecewise linear finite element space with respect to the $H^1(\Omega)$-seminorm, rather than about the FEM itself.

The structure of the paper is as follows. In Section \ref{sec:necessary}, we derive a necessary condition for FEM convergence, i.e. (\ref{sec:intro:est_alpha}) to hold. The condition roughly states that certain structures $\B\subset\Th$, \emph{bands} of neighboring elements forming an alternating pattern such as in Figure \ref{fig:Band} cannot be too long if the maximum angles of elements in $\B$ go to $\pi$ sufficiently fast w.r.t. $h\to 0$. The necessary condition is based on an estimate from below of the error $|u-U|_{H^1(\B)}$ on the band $\B$ -- Theorem \ref{th:H1B_error} for a single band and Theorem \ref{th:mb:H1B_error} for $\Th$ containing multiple bands. Corollaries \ref{cor:H1B_error} and \ref{cor:mb:H1B_error} then state (very roughly, cf. Remark \ref{rem:nec:angles}) that if (\ref{sec:intro:est_alpha}) holds  and if $\B$ has length $L\ge Ch^{2\alpha/5}$ then $\pi-\alpha_K\ge Ch^{3-2\alpha} L$ for all $K\in \B$. This is the case of one band $\B\subset\Th$, for multiple bands this can be improved to $\pi-\alpha_K\ge Ch^{1-\alpha} L$. These results enable us to construct many simple counterexamples to $O(h^\alpha)$ convergence, which is done in Sections \ref{subsec:necessary_examples} and \ref{subsec:necessary_examples:mb}. We note that the Babu\v ska-Aziz counterexample (where $\Th$ consists only of the considered bands, cf. Figure \ref{fig:Babuska-Aziz}) can be obtained by this technique. In fact, we recover the optimal results of \cite{Oswald}. However unlike the Babu\v ska-Aziz counterexample, the regular periodic structure of $\Th$ is not necessary in the presented analysis.

In Section \ref{sec:sufficient} we deal with sufficient conditions for (\ref{sec:intro:est_alpha}) to hold. For simplicity we deal with the $O(h)$ case, i.e. $\alpha=1$ and the general case is then a simple extension (Remark \ref{rem:O_h_alpha}). We split $\Th$ into two parts: $\Th^{1}$ consisting of elements $K$ satisfying (\ref{sec:intro:max}) for a chosen $\alpha_0$ and $\Th^{2}$ consisting of $K$ violating this condition (`degenerating' elements). While on $\Th^{1}$ it is safe to use Lagrange interpolation, on $\Th^{2}$ we use a modified Lagrange interpolation operator to construct a linear function $v_h^K$ on $K$ which is $O(h_K)$-close to $u$ in the $H^1(K)$-seminorm independently of the shape of $K$. The price paid is that $v_h^K$ no longer interpolates $u$ in $A_K$, the maximal-angle vertex of $K$, exactly  but with a small perturbation of the order $O(h_K^2)$, cf. Lemma \ref{lem:th2est}. The question then arises how to connect these piecewise linear functions continuously, if they no longer interpolate the continuous function $u$ exactly in some vertices. Instead of changing the interpolation procedure, we change the interpolated function $u$ itself, so that the Lagrange interpolation of the new function $\tilde u$ corresponds to the modified Lagrange interpolation of $u$ on $\Th^{2}$. Then it is possible to prove that $|u-\Pi_h\tilde u|_1\leq Ch$, hence $|u-U|_1\leq Ch$ by C\' ea's lemma. For this purpose, we introduce the concept of \emph{correction functions} and Sections \ref{subsec:constr_w}, \ref{subsec:clusters} are devoted to constructing the correction functions corresponding to different situations. Roughly speaking, $\Th$ can contain degenerating elements which can form arbitrarily large structures, chains, cf. Figure \ref{fig:chains}, as long as their maximal-angle vertices are not too close to other vertices from $\Th^2$. In general, if the degenerating elements form nontrivial clusters, such as the bands of Section \ref{sec:necessary}, then these can have diameter up to $O(h^{\alpha/2})$. These results are contained in Theorems \ref{thm:suff_main} and \ref{thm:suff_main_cluster}.

We note that while a necessary \emph{and} sufficient condition for any kind of FEM convergence remains an open question, the gap between the derived conditions is small in some special cases, cf. Remark \ref{rem:necandsuff}. For example, in the case of $\Th$ containing one band of length $L$ consisting of sufficiently degenerating elements, the necessary condition for $O(h^\alpha)$ convergence is $L\leq Ch^{2\alpha/5}=Ch^{0.4\alpha}$, while the sufficient condition is $L\leq Ch^{\alpha/2}=Ch^{0.5\alpha}$.

\subsection{Problem formulation and notation}
Let $\Omega\subset\IR^2$ be a bounded polygonal domain with Lipschitz continuous boundary. We treat the following problem: Find $u:\Omega\subset\IR^2\to \IR$ such that
\begin{equation}
-\Delta u=f,\quad u|_{\partial\Omega}=u_D.
\label{cont_prob}
\end{equation}
There are several possibilities how to treat nonhomogeneous Dirichlet boundary conditions in the finite element method, here we consider the standard lifting technique, cf. \cite{Ciarlet}. We note that this choice is not necessary nor important in this paper, as we are essentially interested in the best $H^1$-approximation of $u$ in the discrete space, independent of the specific form of the weak formulation.

Let $g\in H^1(\Omega)$ be the \emph{Dirichlet lift}, i.e. an arbitrary function such that $g|_{\partial\Omega}=u_D$, and write $u=u_0+g$, $u_0\in H^1_0(\Omega)$. Defining $V=H^1_0(\Omega)$ and the bilinear form $a(u,v)=\int_\Omega\nabla u\cdot\nabla v\, \mathrm{d}x$, the corresponding weak form of (\ref{cont_prob}) reads: Find $u_0\in V$ such that
\begin{equation}
a(u+g,v)=(f,v),\quad \forall v\in V.
\nonumber
\end{equation}
The finite element method constructs a sequence of spaces $\{X_h\}_{h\in(0,h_0)}$ on conforming triangulations $\{\Th\}_{h\in(0,h_0)}$ of $\Omega$, where $X_h\subset H^1(\Omega)$ consists of globally continuous piecewise linear functions on $\Th$. Furthermore, let $V_h=X_h\cap V$. We do not assume any properties of $\Th$, since it is the goal of this paper to derive necessary and also sufficient geometric properties of $\Th$ for finite element convergence. 

We choose some approximation $g_h\in X_h$ of $g$, e.g. a piecewise linear Lagrange interpolation of $g$ on $\partial\Omega$ such that $g(x)=0$ in all interior vertices. We seek the FEM solution $U\in X_h$ in the form $U=U_0+g_h$, where $U_0\in V_h$.
The FEM formulation then reads: Find $U_0\in V_h$ such that
\begin{equation}
a(U_0+g_h,v_h)=(f,v_h),\quad \forall v_h\in V_h.
\nonumber
\end{equation}
Formally, we should use the notation $U_h$ instead of $U$ to indicate that $U_h\in X_h$, however we drop the subscript $h$ for simplicity of notation. One must have this in mind e.g. when dealing with estimates of the form $|u-U|_\HO\leq Ch.$

For simplicity of notation, Sobolev norms and seminorms on $\Omega$ will be denoted using simplified notation, e.g. $|\cdot|_1:=|\cdot|_\HO,\ |\cdot|_{2,\infty}:=|\cdot|_{W^{2,\infty}(\Omega)}$, etc. Throughout the paper, $C$ will denote a generic constant independent of $h$ and other geometric quantities of  $\Th$. The constant $C$ has in general different values in different parts of the paper, even within the same chain of inequalities.

\section{A necessary condition for FEM convergence}
\label{sec:necessary}
Throughout Section \ref{sec:necessary} we will assume the following bound on the error. From this error estimate we will derive necessary conditions on the geometry of the triangulations $\Th$ for this estimate to hold.

\medskip
\noindent{\bf Error assumption:} There exists $\alpha\in[0,1]$ independent of $h$ such that for all $u\in H^2(\Omega), \ h\in(0,h_0)$
\begin{equation}
|u-U|_1\leq C(u)h^\alpha.
\label{error_estimate}
\end{equation}

\begin{remark}
\label{rem:err_ass}
Especially of interest are the cases $\alpha=1$ and $\alpha=0$ corresponding to $O(h)$ convergence and (possibly) nonconvergence of the FEM. We note that we can just as well consider the more general case of $|u-U|_\HO\leq C(u)f(h)$ for some (nonincreasing) function $f$. For simplicity we consider only (\ref{error_estimate}) since this coincides well with the concept of orders of convergence.
\end{remark}

\begin{remark}
Assuming (\ref{error_estimate}) holds for all $u\in H^2(\Omega)$ is not necessary, however it is typical for FEM analysis. In Section \ref{sec:necessary}, we will only use (\ref{error_estimate}) for one specific quadratic function (\ref{u_def_quadfunc}). Hence $H^2(\Omega)$ can be replaced by any space that contains quadratic functions.
\end{remark}

\begin{remark}
\label{rem:approx}
We note that $U$ in (\ref{error_estimate}) need not be the FEM solution, since we do not use the FEM formulation in any way in the following. We are essentially dealing only with approximation properties of the space $V_h$: If there exists $U\in V_h$ such that (\ref{error_estimate}) holds, what are the necessary properties of $\Th$.
\end{remark}

\begin{figure}[t]
\begin{center}
\includegraphics[scale=0.7,clip]{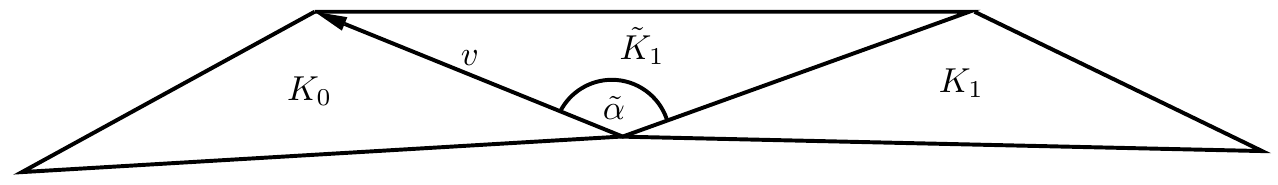}
\caption{Three neighboring elements $K_0,\tilde{K}_1,K_1$.}
\label{fig:3_elements0}
\end{center}
\end{figure}

The basic result upon which we will build the estimates of Section \ref{sec:necessary} is a simple geometric identity concerning gradients of a continuous piecewise linear function $U$ (not necessarily the FEM solution) on a triplet of neighboring elements $K_0,\tilde{K}_1,K_1$ in the configuration from Figure \ref{fig:3_elements0}. We denote $U_0:=U|_{K_0}$, $U_1:=U|_{K_1}$ and $\tilde U_1:=U|_{\tilde K_1}$. The key observation is that if $U_0, U_1$ are given, then $\tilde U_1$ is uniquely determined due to inter-element continuity. We note that the following result is not an estimate, it is an equality, hence optimal.

\begin{lemma}
\label{lem:angles}
Let $\xi$ be the angle between the vectors $\grad(U_0-U_1)$ and $v$, where $v$ corresponds to the common edge of $K_0,\tilde K_1$, cf. Figure \ref{fig:3_elements0}. Then
\begin{equation}
|\grad(\tilde{U}_1-U_1)|=\frac{\cos(\xi)}{\sin(\pi-\talpha)}|\grad(U_0-U_1)|,
\label{lem:angles_eq}
\end{equation}
where $\tilde{\alpha}$ is the maximum angle of $\tilde K_1$.
\end{lemma}
\proof
\begin{figure}[t]
\begin{center}
\includegraphics[scale=0.7,clip]{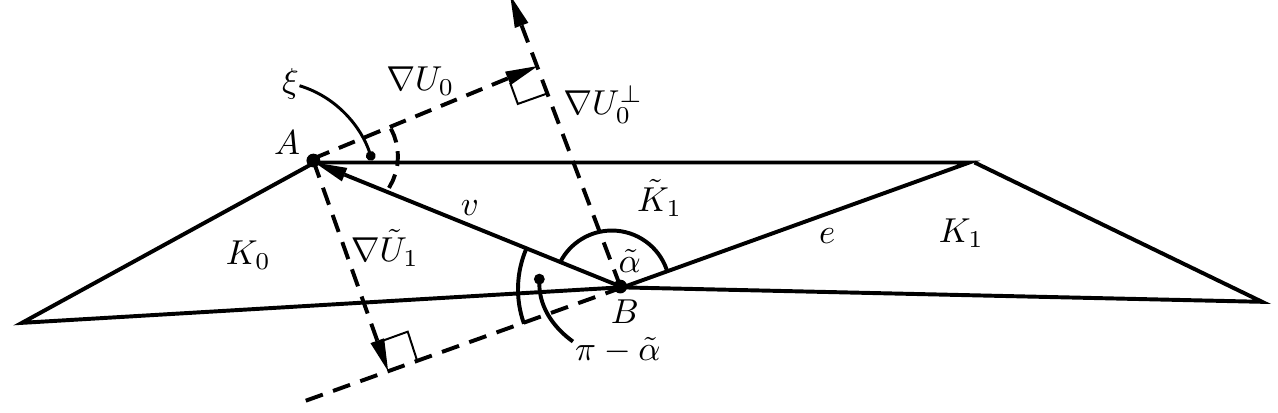}
\caption{Proof of Lemma \ref{lem:angles}.}
\label{fig:3_elements}
\end{center}
\end{figure}
By globally subtracting the function $U_1$ (i.e. its extension to the whole $\IR^2$) from all functions on $K_1,\tilde{K}_1, K_0$, we can assume that $U_1\equiv 0$ without loss of generality. From the continuity of $U$, we have $U_0(A)=\tilde{U}_1(A)$, where $A$ is the common vertex of $K_0$ and $\tilde K_1$ denoted in Figure \ref{fig:3_elements}. Furthermore, by continuity and since $U_1\equiv 0$, then $\tilde U_1=0$ on the edge $e$. Therefore 
\begin{equation}
|\grad(\tilde{U}_1-U_1)|=|\grad\tilde{U}_1|=\frac{|\tilde U_1(A)|}{|v|\sin(\pi-\talpha)}.
\label{lem:angles1}
\end{equation} 
On the other hand, by continuity and since $U_1\equiv 0$, we have $U_0(B)=0$ and $U_0=0$ on the line $l$ passing through $B$ in the direction $\grad U_0^\perp=\grad (U_0-U_1)^\perp$. Therefore,
\begin{equation}
|\grad(U_0-U_1)|=|\grad U_0|=\frac{|\tilde U_1(A)|}{|v|\cos(\xi)}.
\label{lem:angles2}
\end{equation} 
Expressing $|\tilde U_1(A)|$ from (\ref{lem:angles2}) and substituting into (\ref{lem:angles1}), we obtain the desired result.
\qed

\begin{remark}
\label{rem:angles}
Since $(a,b)=|a||b|\cos\alpha$, where $\alpha$ is the angle between vectors $a,b$, we can reformulate (\ref{lem:angles_eq}) as 
\begin{equation}
|\grad(\tilde{U}_1-U_1)| =\frac{1}{\sin(\pi-\talpha)} \big|\big(\tfrac{1}{|v|}v,\grad(U_0-U_1)\big)\big|.
\label{rem:angles_eq}
\end{equation}
In other words, on the right-hand side of (\ref{rem:angles_eq}), we have the magnitude of the projection of $\grad(U_0-U_1)$ into the direction given by the unit vector $\tfrac{1}{|v|}v$.  
\end{remark}

\medskip
Lemma \ref{lem:angles} states that if $|\grad(U_0-U_1)|$ is nonzero, then $|\grad(\tilde{U}_1-U_1)|$ is huge, since it is magnified by the factor $\tfrac{1}{\sin(\pi-\talpha)}$ which tends to $\infty$ as $\talpha\to \pi$. If $u$ is e.g a reasonable quadratic function and if $\grad U_1$ is a good approximation of $\grad u|_{K_1}$ then  $|\grad U_1|\approx 1$, hence $\grad\tilde{U}_1$ will be huge, hence a bad approximation of $\grad u|_{\tilde K_1}$, which is reasonably bounded. 

\begin{figure}[t]
\begin{center}
\includegraphics[scale=0.7,clip]{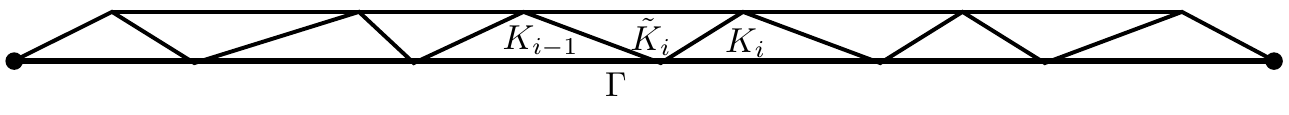}
\caption{A band of elements $\B$ with edge $\Gamma$.}
\label{fig:Band}
\end{center}
\end{figure}

We will proceed as follows: We will choose an exact solution $u$ with nonzero second derivatives (a quadratic function). If $|u-U|_1$ is small, then we can expect $|\grad(U_0-U_1)|\neq 0$, since this is essentially an  approximation of second order derivatives. Therefore $|\grad(\tilde U_1-U_1)|$ will be enormous on degenerating triangles due to the factor $1/\sin(\pi-\tilde \alpha)$. Hence $\grad\tilde U_1$ cannot be a good approximation of $u|_{\tilde K_1}$, therefore $|u-U|_1$ cannot be small. From this we can obtain restrictions e.g. on $\sin(\pi-\tilde \alpha)$ for (\ref{error_estimate}) to hold. The main difficulty is that the assumption of $|u-U|_1$ being small does not imply $|\grad(U_0-U_1)|\neq 0$ for any given triplet of elements such as in Figure \ref{fig:3_elements0}, we know this only ``on average", in the $L^2(\Omega)$-sense due to the $H^1(\Omega)$ error estimate. However, if we connect the considered triplets of elements into larger structures (bands $\B$, cf. Figure \ref{fig:Band}), then if these bands are large enough, for some triplet of elements in $\B$ necessarily $|\grad(U_0-U_1)|\neq 0$ and we can apply Lemma \ref{lem:angles}.

Our goal will be to prove that if (\ref{error_estimate}) holds and $\B$ is long enough, then the differences $|\grad(U_{i+1}-U_i)|$ must be nonzero on average, hence $\grad\tilde U_i$ is a bad approximation of $\grad u|_{\tilde K_i}$ and therefore the global error cannot be $O(h^\alpha)$ on the elements $\tilde{K}_i$. Such contradictions allow us to formulate a necessary condition for (\ref{error_estimate}) and to construct counterexamples to $O(h)$ and other convergences.

\medskip
\begin{definition}
\label{def:band}
We define a band $\B\subset\Th$ as the union $\B=(\cup_{i=0}^N K_i)\cup(\cup_{i=1}^N \tilde K_i)$ such that $K_{i-1},\tilde K_i,K_i$ for $i=1,\ldots,N$ are neighboring elements forming a triplet as in Lemma \ref{lem:angles}, such that their maximum angles form an alternating pattern as in Figure \ref{fig:Band}. We denote the union of longest edges of all $K_i$ by $\Gamma$.
\end{definition}

\begin{remark}
\label{rem:straightB}
In the following, we assume that the edge $\Gamma$ of band $\B$ lies on a line, i.e. $\B$ is `straight'. This makes the notation and proofs less technical, however it is not a necessary assumption. We will indicate in the relevant places how the proofs should be modified to allow for $\B$ `curved'.
\end{remark}

\begin{lemma}
\label{lem:L2G}
Let $\B$ be a band as in definition $\ref{def:band}$ and let the error estimate (\ref{error_estimate}) hold. Then for sufficiently small $h$
\begin{equation}
\|u-U\|_\LG\leq Ch^\alpha,
\label{lem:L2G_est}
\end{equation}
where $C$ depends only on $\Omega$ and the constant in (\ref{error_estimate}).
\end{lemma}
\proof
The proof is essentially a straightforward application of the trace inequality. However, we must be careful, since in the standard trace inequality the constant depends heavily on the geometry of the domain. Therefore, we give a detailed proof of (\ref{lem:L2G_est}). We proceed similarly as in \cite{Grisvard}.

\begin{figure}[t]
\begin{center}
\includegraphics[scale=0.6,clip]{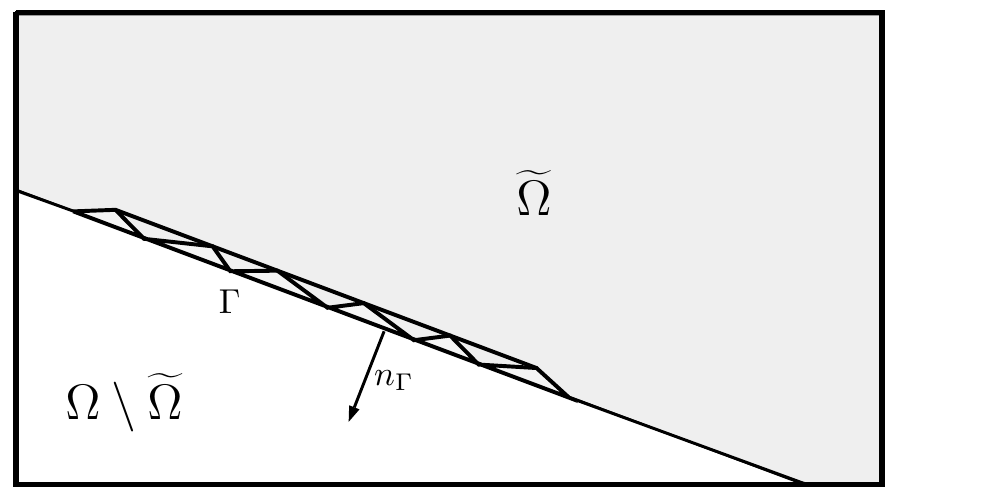}
\caption{Proof of Lemma \ref{lem:L2G} -- splitting of $\Omega$ by $\Gamma$.}
\label{fig:Trace_ineq}
\end{center}
\end{figure}

For simplicity, we denote $e:=u-U$. We define the constant vector function $\bmu=-\n_\Gamma$, where $\n_\Gamma$ is the unit outer normal to $\B$ on $\Gamma$, cf. Figure \ref{fig:Trace_ineq}. Then $\bmu$ is the unit outer normal to $\widetilde{\Omega}$ on $\Gamma$, where $\widetilde{\Omega}\subset\Omega$ is the subdomain containing $\B$ defined by the line on which $\Gamma$ lies. On one hand, we have
\begin{equation}
\int_{\widetilde{\Omega}} \grad(e^2)\cdot\bmu\dx=2\int_{\widetilde{\Omega}} e\grad e\cdot\bmu\dx.
\nonumber
\end{equation}
On the other hand, by Green's theorem,
\begin{equation}
\int_{\widetilde{\Omega}} \grad (e^2)\cdot\bmu\dx=\int_{\pd\widetilde{\Omega}} e^2\bmu\cdot\n\dS,
\nonumber
\end{equation}
since $\diver\bmu=0$. As $\bmu\cdot\n=1$ on $\Gamma$, by combining these two equalities we get
\begin{equation}
\int_{\Gamma} e^2\dS\leq 2\int_{\widetilde{\Omega}} e\grad e\cdot\bmu\dx-\int_{\pd\widetilde{\Omega}\cap\pd\Omega} e^2\bmu\cdot\n\dS,
\label{lem:L2G_3}
\end{equation}
where we have omitted the boundary integral over $\pd\widetilde{\Omega}\setminus(\pd\Omega\cup\Gamma)$, which is nonnegative since $\bmu\cdot\n=1$ on this part of $\pd\widetilde{\Omega}$. The second right-hand side term in (\ref{lem:L2G_3}) can be estimated as
\begin{equation}
\int_{\pd\widetilde{\Omega}\cap\pd\Omega} e^2\bmu\cdot\n\dS\leq \int_{\pd\Omega} e^2\dS\leq Ch^4,
\label{lem:L2G_4}
\end{equation}
since on $\pd\Omega$, $U$ is the piecewise linear Lagrange interpolation of $u$, hence $|e|=|u-U|\leq Ch^2$ on $\partial\Omega$. The first right-hand side term of (\ref{lem:L2G_3}) can be estimated using Young's inequality by
\begin{equation}
2\int_{\widetilde{\Omega}} e\grad e\cdot\bmu\dx\leq\int_{\Omega} e^2+|\grad e\cdot\bmu|^2\dx\leq \|e\|^2_\LO+|e|^2_\HO.
\label{lem:L2G_5}
\end{equation}
Combining (\ref{lem:L2G_3})--(\ref{lem:L2G_5}) with (\ref{error_estimate}), we get
\begin{equation}
\|e\|^2_\LG\leq \|e\|^2_\LO+Ch^4+Ch^{2\alpha}.
\label{lem:L2G_6}
\end{equation}

In (\ref{lem:L2G_6}), it remains to estimate $\|e\|^2_\LO$ by $|e|^2_\HO$. Ordinarily, one could simply use Poincar\'{e}'s inequality, however $e\notin H^1_0(\Omega)$. Nonetheless, as mentioned above, $e|_{\pd\Omega}$ is of the order $O(h^2)$, which allows us to obtain a similar estimate to Poincar\'{e}'s inequality. Using Green's theorem and the trivial identity $\frac{\pd x_1}{\pd x_1}=1$, we estimate
\begin{equation}
\|e\|^2_\LO=\int_\Omega 1e^2\dx=\int_{\pd\Omega}x_1 n_1 e^2\dS-\int_\Omega x_12e\frac{\pd e}{\pd x_1}\dx.
\label{lem:L2G_7}
\end{equation}
Without loss of generality, let $0\in\Omega$. Therefore we can estimate (\ref{lem:L2G_7}) using Young's inequality as
\begin{equation}
\begin{split}
\|e\|^2_\LO&\leq\diam\Omega\|e\|^2_\LG+2\diam\Omega\|e\|_\LO|e|_\HO\\ &\leq\diam\Omega\|e\|^2_\LG+\frac{1}{2}\|e\|^2_\LO+(\diam\Omega)^2|e|^2_\HO.
\nonumber
\end{split}
\end{equation}
Therefore,
\begin{equation}
\|e\|^2_\LO\leq 2\diam\Omega\|e\|^2_\LG+2(\diam\Omega)^2|e|^2_\HO\leq Ch^4+Ch^{2\alpha}.
\label{lem:L2G_9}
\end{equation}
Combining (\ref{lem:L2G_6}) and (\ref{lem:L2G_9}) gives us estimate (\ref{lem:L2G_est}) for sufficiently small $h$.
\qed

\begin{remark}
\label{rem:curvedB_1}
If we do not assume that $\Gamma$ lies on a straight line, cf. Remark \ref{rem:straightB}, we need to find a constant vector function $\mu$ such that $\mu\cdot n\ge\mu_0$ for some $\mu_0>0$. This can be done e.g. by taking $\mu=n_{K_i}^\Gamma$ for some $i\in\{0,\ldots,N\}$, i.e. the normal to $K_i$ on $\Gamma$, and assuming that $n_{K_i}^\Gamma\cdot n_{K_j}^\Gamma\ge\mu_0>0$ for all $j\neq i$, i.e. $\B$ does not `bend too much'. The constant $\mu_0$ then figures in the left-hand sides of (\ref{lem:L2G_3}) and (\ref{lem:L2G_6}). Consequently, we get the constant $1/\mu_0$ in the right-hand side of (\ref{lem:L2G_est}).
\end{remark}

\medskip

In the following, we shall assume that the exact solution of (\ref{cont_prob}) is 
\begin{equation}
u(x_1,x_2)=x_1^2+x_2^2,
\label{u_def_quadfunc}
\end{equation} 
i.e. $f=-4$ with corresponding Dirichlet boundary conditions. If $L$ is the length of $\Gamma$, then $\Gamma$ can be parametrized by $t\in(0,L)$ as
\begin{equation}
\Gamma=\{x=a+tg, t\in(0,L)\},
\nonumber
\end{equation} 
where $a$ is an endpoint of $\Gamma$ and $|g|=1$ is a vector in the direction of $\Gamma$. The restriction $u|_\Gamma$ can then also be parametrized by $t\in(0,L)$. We denote the corresponding function $u_\Gamma$ and observe that 
\begin{equation}
u_\Gamma(t)=u(a+tg)=t^2+p_1(t), 
\label{eq:ug}
\end{equation}
where $p_1\in P^1{(0,L)}$, i.e. a linear function. We have the following approximation result:

\begin{lemma}
\label{lem:ug_proj}
Let $u_\Gamma\in P^2{(0,L)}$ be given by (\ref{eq:ug}). Let $\Pi^1_{(0,L)}$ be the $L^2(0,L)$-orthogonal projection onto $P^1(0,L)$. Then
\begin{equation}
\|u_\Gamma-\Pi^1_{(0,L)}u_\Gamma\|_{L^2(0,L)}=\tfrac{1}{6\sqrt{5}}L^{5/2}.
\label{lem:ug_proj_eq}
\end{equation}
\end{lemma}
\proof
For convenience, $u_\Gamma$ can be rewritten in the form $u_\Gamma(t)=(t-L/2)^2+\tilde{p}_1(t)$, where $\tilde{p}_1\in P^1(0,L)$. Then $u_\Gamma-\Pi^1_{(0,L)}u_\Gamma=(t-L/2)^2-\Pi^1_{(0,L)}(t-L/2)^2$, since $\Pi^1_{(0,L)}\tilde{p}_1=\tilde{p}_1$. 

On $(0,L)$, the quadratic function $(t-L/2)^2$ is symmetric with respect to $L/2$, hence its projection must be a constant function. It can be therefore easily computed that $\Pi^1_{(0,L)}(t-L/2)^2=\tfrac{1}{12}L^3$. Therefore the norm in (\ref{lem:ug_proj_eq}) can be straightforwardly computed, giving the desired result. To save space, we omit the elementary, yet lengthy calculations.
\qed

\begin{corollary}
\label{cor:ug_proj}
Let $\alpha\in[0,1]$ and $L\ge Ch^{2\alpha/5}$. Then
\begin{equation}
\|u_\Gamma-\Pi^1_{(0,L)}u_\Gamma\|_{L^2(0,L)}\ge Ch^\alpha.
\label{cor:ug_proj_eq}
\end{equation}
\end{corollary}

\begin{remark}
\label{rem:ug_proj}
Specifically of interest are the cases $\alpha=1$, corresponding to $L\ge Ch^{2/5}$ and $\alpha=1$, corresponding to $L\ge C$. These two cases will lead to necessary conditions and counterexamples to $O(h)$ convergence, and convergence of the FEM, respectively.
\end{remark}

\begin{remark}
\label{rem:curvedGamma2}
If $\Gamma$ does not lie on a straight line, cf. Remark \ref{rem:straightB}, the length-parametrized restriction $u_\Gamma:(0,L)\to\IR$ of $u|_\Gamma$ will be continuous, piecewise quadratic. However, if $\Gamma$ is close to a straight line, then $u_\Gamma$ will be close to the globally quadratic function $t^2+p_1(t)$ of (\ref{eq:ug}). We can then estimate
\begin{equation}
\begin{split}
\|&u_\Gamma-\Pi^1_{(0,L)}u_\Gamma\|_{L^2(0,L)}\\
&\ge \bigl|\|t^2+p_1(t)-\Pi^1_{(0,L)}u_\Gamma\|_{L^2(0,L)} -\|u_\Gamma-t^2+p_1(t)\|_{L^2(0,L)}\bigr|.
\label{rem:curvedGamma2:1}
\end{split}
\end{equation}
For the first right-hand side term we have
\begin{equation}
\|t^2+p_1(t)-\Pi^1_{(0,L)}u_\Gamma\|_{L^2(0,L)} \ge\|t^2+p_1(t)-\Pi^1_{(0,L)}(t^2+p_1(t))\|_{L^2(0,L)} =\tfrac{1}{6\sqrt{5}}L^{5/2}
\nonumber
\end{equation}
by Lemma \ref{lem:ug_proj}. If we assume e.g. that 
\begin{equation}
\|u_\Gamma-t^2+p_1(t)\|_{L^2(0,L)}\leq \tfrac{1}{12\sqrt{5}}L^{5/2}
\label{rem:curvedGamma2:2}
\end{equation}
then by (\ref{rem:curvedGamma2:1}) we get the statement of Lemma \ref{lem:ug_proj} as an estimate from below with the lower bound $\tfrac{1}{12\sqrt{5}}L^{5/2}$. Condition (\ref{rem:curvedGamma2:2}) is effectively a condition on how $\Gamma$ deviates from a straight line.
\end{remark}

Now we focus on the discrete solution on $\Gamma$. Again, if $\Gamma$ lies on a straight line, then the length-parametrized restriction $U_\Gamma:=U|_\Gamma\in L^2(0,L)$ is a continuous piecewise linear function on the partition of $(0,L),$ i.e. $\Gamma$, induced by $\Th$. Hence $U_\Gamma\in H^1(0,L)$.

For convenience, to deal with functions from $H^1_0(0,L)$ instead of $H^1(0,L)$, we define $\bar U\in P^1(0,L)$ such that $\bar U(0)=U_\Gamma(0)$ and $\bar U(L)=U_\Gamma(L)$. Let
\begin{equation}
w:=\UG-\bar U,
\label{w_def}
\end{equation}
then $w(0)=w(L)=0$. We have the following estimate.

\begin{lemma}
\label{lem:wH1est}
Let $L\ge C_Lh^{2\alpha/5}$, where $C_L$ is sufficiently large with respect to the constant $C(u)$ from (\ref{error_estimate}). Then
\begin{equation}
|w|_\HoL\ge C L^{3/2}.
\label{lem:wH1est_est}
\end{equation}
\end{lemma}
\proof
By the reverse triangle inequality, Lemma \ref{lem:ug_proj} and \ref{lem:L2G},
\begin{equation}
\begin{split}
\|w\|_\LoL&\ge\big| \|\uG-\bar U\|_\LoL-\|\uG-\UG\|_\LoL \big|\\
&\ge \big| \|\uG-\Pi^1_{(0,L)}\uG\|_\LoL-\|\uG-\UG\|_\LoL \big|\\ &\ge |\tfrac{1}{6\sqrt{5}}L^{5/2}-Ch^\alpha|\ge \tfrac{1}{12\sqrt{5}} L^{5/2},
\label{lem:wH1est_1}
\end{split}
\end{equation}
since $L^{5/2}\ge C_L^{5/2}h^\alpha$ and this term dominates the error term $u-\UG$ for $C_L$ sufficiently large w.r.t. the constant in (\ref{lem:L2G_est}), i.e. (\ref{error_estimate}).

To obtain (\ref{lem:wH1est_est}), we use Poincar\'e's inequality on $(0,L)$. This is possible, since $w$ is continuous, piecewise linear and $w(0)=w(L)=0$, hence $w\in H^1_0(0,L)$. Since $L$ is not fixed, one must use the appropriate scaling from the Poincar\'e inequality on the fixed interval $(0,1)$. From (\ref{lem:wH1est_1}) we obtain
\begin{equation}
\tfrac{C}{2}L^{5/2}\leq \|w\|_\LoL\leq C_PL|w|_\HoL, 
\nonumber
\end{equation}
where $C_P$ is the constant from Poincar\'e's inequality on $(0,1)$. Cancelling $L$ from both sides gives us (\ref{lem:wH1est_est}).
\qed

\medskip
Since $w$ is piecewise linear on $(0,L)$, $w^\prime$ is piecewise constant on the corresponding partition  $(0,L)=\overline{\cup_{i=0}^N I_i}$ induced on $\Gamma$ by $\Th$. We denote $h_i=|I_i|$ and $\wp_i=(w|_{I_i})^\prime$. Then Lemma \ref{lem:wH1est} can be rewritten as
\begin{equation}
\sum_{i=0}^N h_i{\wp_i}^2\ge CL^3.
\label{wineq1}
\end{equation}
Furthermore, since $w(0)=w(L)=0$, we have $\int_0^L \wp(x)\dx=0$, i.e.
\begin{equation}
\sum_{i=0}^N h_i{\wp_i}=0.
\label{wineq2}
\end{equation}
From these two inequalities, we derive an estimate of $\wp_i-\wp_{i-1}$, i.e. how much $w$ differs from a globally linear function. This is already related to the right-hand side of (\ref{lem:angles_eq}).

\begin{lemma}
\label{lem:differences}
Let $\{h_i\}_{i=0}^N,\{a_i\}_{i=0}^N\in\IR^N$ be such that $h_i>0,\ i=0,\ldots,N$ and
\begin{equation}
\sum_{i=0}^N h_i a_i^2=A,\quad \sum_{i=0}^N h_i a_i=0,\quad \sum_{i=0}^N h_i=L,
\label{lem:differences_assumpt}
\end{equation} 
where $A\ge 0$. Then
\begin{equation}
\sum_{i=1}^N (a_i-a_{i-1})^2\ge\frac{A}{LN}.
\label{lem:differences_est}
\end{equation}
\end{lemma}
\proof
Assume on the contrary, that there exists $\{a_i\}_{i=0}^N$ such that (\ref{lem:differences_assumpt}) holds and
\begin{equation}
\sum_{i=1}^N (a_i-a_{i-1})^2<\frac{A}{LN}.
\label{lem:differences_1}
\end{equation}
Let $i_0$ be such that $|a_{i_0}|=\max_{i=0,\cdots,N} |a_i|$. Without loss of generality assume e.g. $a_{i_0}>0$. Since all $h_i>0$ and $\sum_{i=0}^n h_ia_i=0$, there exists $i_1$ such that $a_{i_1}<0$. Without loss of generality assume e.g. $i_0<i_1$. Then for all $i$
\begin{equation}
\begin{split}
|&a_i|\leq a_{i_0}< a_{i_0}-a_{i_1}=(a_{i_0}-a_{i_0+1})+\ldots+(a_{i_1-1}-a_{i_1})\\
&\leq \sum_{k=1}^N|a_k-a_{k-1}|\leq N^{1/2}\Big(\sum_{k=1}^N|a_k-a_{k-1}|^2\Big)^{1/2}<A^{1/2}L^{-1/2}
\end{split}
\label{lem:differences_2}
\end{equation}
due to the Cauchy-Schwartz inequality and (\ref{lem:differences_1}). From (\ref{lem:differences_assumpt}), we have $\sum_{i=0}^N h_i a_i^2=A$. On the other hand, using (\ref{lem:differences_2}),
\begin{equation}
A=\sum_{i=0}^N h_i a_i^2< \sum_{i=0}^N h_i (A^{1/2}L^{-1/2})^2=AL^{-1}\sum_{i=0}^N, h_i=A,
\nonumber
\end{equation}
i.e. $A<A$ which is a contradiction.
\qed

\begin{lemma}
\label{lem:wdif}
Let $L\ge C_Lh^{2\alpha/5}$, where $C_L$ is sufficiently large. Then
\begin{equation}
\sum_{i=1}^N (U_i^\prime-U_{i-1}^\prime)^2\ge CL^2N^{-1},
\label{lem:wdif_est}
\end{equation} 
where $U_i^\prime=(U_\Gamma|_{I_i})^\prime$.
\end{lemma}
\proof
We apply Lemma \ref{lem:differences} to inequalities (\ref{wineq1}), (\ref{wineq2}) with $a_i:=\wp_i$. Then $a_i-a_{i-1}=\wp_i-\wp_{i-1}=U_i^\prime-U_{i-1}^\prime$, since $\bar U^\prime$ is constant. The assumption $L\ge C_Lh^{2\alpha/5}$ is needed because (\ref{wineq1}) follows from Lemma \ref{lem:wH1est}.
\qed

\begin{remark}
\label{rem:C_suff_large}
The phrase ``$C_L$ \textit{sufficiently large}" in Lemmas (\ref{lem:wH1est}), (\ref{lem:wdif}) should read in full ``$C_L$ \textit{sufficiently large w.r.t. the constant $C(u)$ from (\ref{error_estimate}) with $u$ given by (\ref{u_def_quadfunc})}". This can be seen from the proof of Lemma \ref{lem:wH1est}, from which it carries on to further results, where we will use the shortened version ``$C_L$ \textit{sufficiently large}". We note that the dependence of the necessary $C_L$ on $C(u)$ is such that when $C(u)\to 0$ the minimal necessary $C_L$ also goes to zero. We use this fact in Counterexample 3.
\end{remark}

\subsection{Main estimate -- one band}
\label{sec:main_est:one_band}
Let $\B=(\cup_{i=0}^N K_i)\cup(\cup_{i=1}^N \tilde K_i)$ be the band with edge $\Gamma$, cf. Figure \ref{fig:Band}. We denote $U_i=U|_{K_i}$ and $\tilde U_i=U|_{\tilde K_i}$. By $g$ we denote the unit vector in the direction of $\Gamma$, i.e. $|g|=1$ and by $\alpha_i,\talpha_i$ we denote the maximum angles of $K_i$ and $\tK_i$, respectively. 

For $U_i^\prime$ from Lemma \ref{lem:wdif} we have $U_i^\prime=\grad U_i\cdot g$, since by definition $U_i^\prime$ is the derivative of $U_i$ in the direction $g$. Hence we can reformulate estimate (\ref{lem:wdif_est}) as
\begin{equation}
\sum_{i=1}^N \big|(\grad U_i-\grad U_{i-1})\cdot g\big|^2\ge CL^2N^{-1}.
\label{lem:wdif_est2}
\end{equation}
 
We denote $\grad u_i=\grad u(\mcC_i)$, where $\mcC_i$ is the centroid of $K_i$ and similarly $\grad \tilde u_i=\grad u(\tilde \mcC_i)$, where $\tilde \mcC_i$ is the centroid of $\tilde K_i$. Trivially, 
\begin{equation}
\begin{split}
|\grad u(x)-\grad u_i|\leq Ch,\quad \forall x\in K_i,\\
|\grad u(x)-\grad \tilde u_i|\leq Ch,\quad \forall  x\in \tilde K_i.
\end{split}
\label{grad_ui_approx}
\end{equation}

Now we will prove the main theorem of this section. We show that if $|u-U|_1\leq Ch^\alpha$, then the error on $\B$ can be split into two parts. The first ``dominant" part $\mcA_1$ can be bounded from below, while the remainder $\mcA_2$ is of the order $O(h^\alpha)$. The estimate for $A_1$ depends on the geometry of $\B$ and can be made arbitrarily large e.g. for the maximal angles in $\B$ going to $\pi$ sufficiently fast. Thus contradictions with the $O(h^\alpha)$ bound of the error can be created, producing counterexamples to $O(h^\alpha)$ convergence as in Section \ref{subsec:necessary_examples}.

\begin{theorem}
\label{th:H1B_error}
Let $u$ be given by (\ref{u_def_quadfunc}) and let $U\in X_h$ satisfy the error estimate (\ref{error_estimate}). Let the band $\B=(\cup_{i=0}^N K_i)\cup(\cup_{i=1}^N \tilde K_i)$ be such that there exist constants $C_K, C_\alpha$ independent of $h$ such that for all $i=1,\ldots,N$
\begin{equation}
\begin{split}
|\tilde K_i|&\leq C_K|K_i|,\\
|\tilde K_i|&\leq C_K|K_{i-1}|,\\
\sin(\pi-\alpha_{i-1})&\leq C_\alpha\sin(\pi-\tilde\alpha_{i}).
\label{th:H1B_error_ass2}
\end{split}
\end{equation}
Finally, let $L\ge C_Lh^{2\alpha/5}$, where $C_L$ is sufficiently large. Then 
\begin{equation}
|u-U|_{H^1(\B)}\ge \mcA_1-\mcA_2,
\label{th:H1B_error_est1}
\end{equation}
where 
\begin{equation}
\begin{split}
\mcA_1&\ge C_1\min_{i=1,\ldots,N}\frac{1}{\sin(\pi-\talpha_i)}\min_{i=1,\ldots,N}|\tilde K_i|^{1/2} LN^{-1/2},\\
|\mcA_2|&\leq C_2h^\alpha.
\end{split}
\label{th:H1B_error_est2}
\end{equation}
\end{theorem}
\proof
We consider the error only on the elements $\tilde K_i$:
\begin{equation}
\begin{split}
|u&-U|_{H^1(\B)}\ge \Big(\sum_{i=1}^N\int_{\tilde K_i}|\grad u-\grad U|^2\dx\Big)^{1/2} =\mcA_1-\mcA_2,
\end{split}
\label{th:H1B_error_eq1}
\end{equation}
where the splitting of the error is given by
\begin{equation}
\begin{split}
\mcA_1&:=\Big(\sum_{i=1}^N \frac{1}{\sin^2(\pi-\talpha_i)}\big|(\grad U_i-\grad U_{i-1})\cdot g\big|^2|\tilde K_i|\Big)^{1/2},\\
\mcA_2&:=\Big(\sum_{i=1}^N\int_{\tilde K_i}|\grad u-\grad U|^2\dx\Big)^{1/2}-\mcA_1.
\label{th:H1B_error_eq2}
\end{split}
\end{equation}
Due to (\ref{lem:wdif_est2}) we immediately have
\begin{equation}
\mcA_1\ge \min_{i=1,\ldots,N}\frac{1}{\sin(\pi-\talpha_i)}\min_{i=1,\ldots,N}|\tilde K_i|^{1/2} \big(CL^2N^{-1}\big)^{1/2},
\label{th:H1B_error_eq3}
\end{equation}
which is the first inequality in (\ref{th:H1B_error_est2}). It remains to estimate $\mcA_2$.

Since $\grad U|_{\tilde{K}_i}=\grad \tilde{U}_i$, a series of triangle inequalities gives us
\begin{equation}
\begin{split}
\mcA_2&\leq \Bigl(\sum_{i=1}^N\int_{\tilde K_i}|\grad u-\grad \tilde u_i|^2\dx\Bigr)^{1/2}\!\!+\Bigl(\sum_{i=1}^N\int_{\tilde K_i}|\grad \tilde u_i-\grad \tilde U_i|^2\dx\Bigr)^{1/2}\!\!\! -\!\mcA_1\\
&=\Bigl(\sum_{i=1}^N\int_{\tilde K_i}|\grad u-\grad \tilde u_i|^2\dx\Bigr)^{1/2} \!\!+\Big(\sum_{i=1}^N|\grad \tilde u_i-\grad \tilde U_i|^2|\tilde K_i|\Big)^{1/2}\!\!-\mcA_1\\
&\leq \Big(\sum_{i=1}^N\int_{\tilde K_i}|\grad u-\grad \tilde u_i|^2\dx\Big)^{1/2}+
\Big(\sum_{i=1}^N|\grad \tilde u_i-\grad u_i|^2|\tilde K_i|\Big)^{1/2}\\ &\quad+\Big(\sum_{i=1}^N|\grad u_i-\grad U_i|^2|\tilde K_i|\Big)^{1/2}\!\! +\Big(\sum_{i=1}^N|\grad U_i-\grad \tilde U_i|^2|\tilde K_i|\Bigr)^{1/2}\!\!-\!\mcA_1
\\ &=:(A)+(B)+(C)+(D)-\mcA_1.
\label{th:H1B_error_eq4}
\end{split}
\end{equation}
We estimate the individual terms $(A)-(D)$ of (\ref{th:H1B_error_eq4}). Due to (\ref{grad_ui_approx}),
\begin{equation}
(A)\leq Ch|\B|^{1/2}\leq Ch|\Omega|^{1/2}\leq Ch\leq Ch^\alpha,
\label{th:H1B_error_eq5}
\end{equation}
since $\alpha\in [0,1]$. Similarly,
\begin{equation}
(B)=\Big(\sum_{i=1}^N\big|\grad u(\tilde \mcC_i)-\grad u(\mcC_i)\big|^2|\tilde K_i|\Big)^{1/2}\leq Ch|\B|^{1/2}\leq Ch\leq Ch^\alpha,
\label{th:H1B_error_eq6}
\end{equation}
since $\grad u$ is a linear function due to (\ref{u_def_quadfunc}) and the centroids of the neighboring elements $K_i,\tilde K_i$ satisfy $|\mcC_i-\tilde \mcC_i|\leq 2h$.

Under the assumption $|\tilde K_i|\leq C_K|K_i|$ for all $i$, we have
\begin{align}
(C)&\leq\! \Bigl(C_K\sum_{i=1}^N|\grad u_i-\grad U_i|^2|K_i|\Bigr)^{1/2} \!\!=\!\Bigl(C_K\sum_{i=1}^N\int_{K_i}|\grad u_i-\grad U|^2\dx\Bigr)^{1/2}\nonumber\\
&\leq\! \Bigl(C_K\!\sum_{i=1}^N\int_{K_i}\!|\grad u_i-\grad u|^2\dx\Bigr)^{1/2} \!\!+\!\Bigl(C_K\!\sum_{i=1}^N\int_{K_i}\!|\grad u-\grad U|^2\dx\Bigr)^{1/2}\nonumber\\
&\leq Ch|\B|^{1/2}+C_K^{1/2}|u-U|_1\leq Ch^\alpha
\label{th:H1B_error_eq7}
\end{align}
due to (\ref{grad_ui_approx}) and (\ref{error_estimate}).

Using Lemma \ref{lem:angles} (or Remark \ref{rem:angles}) and denoting by $v_{i-1}$ the unit vector corresponding to the common edge of $K_{i-1},\tilde{K}_i$, we can rewrite $(D)$ as
\begin{equation}
\begin{split}
(D)&=\Big(\sum_{i=1}^N\frac{1}{\sin^2(\pi-\talpha_i)} \big|v_{i-1}\cdot\grad(U_i-U_{i-1})\big|^2|\tilde K_i|\Big)^{1/2}\\
&\leq \Big(\sum_{i=1}^N \frac{1}{\sin^2(\pi-\talpha_i)}\big|(\grad U_i-\grad U_{i-1})\cdot g\big|^2|\tilde K_i|\Big)^{1/2}\\
&\quad +\Big(\sum_{i=1}^N \frac{1}{\sin^2(\pi-\talpha_i)}\big|(\grad U_i-\grad U_{i-1})\cdot (v_{i-1}-g)\big|^2|\tilde K_i|\Big)^{1/2}\\
&=:\mcA_1+(E).
\label{th:H1B_error_eq8}
\end{split}
\end{equation}
We notice that we have obtained the term $\mcA_1$ in the estimate of $(D)$, cf. (\ref{th:H1B_error_eq2}). This will cancel with the final term $-\mcA_1$ in (\ref{th:H1B_error_eq4}), while all the remaining terms are estimated by $Ch^\alpha$.

\begin{figure}[t]
\begin{center}
\includegraphics[scale=0.7,clip]{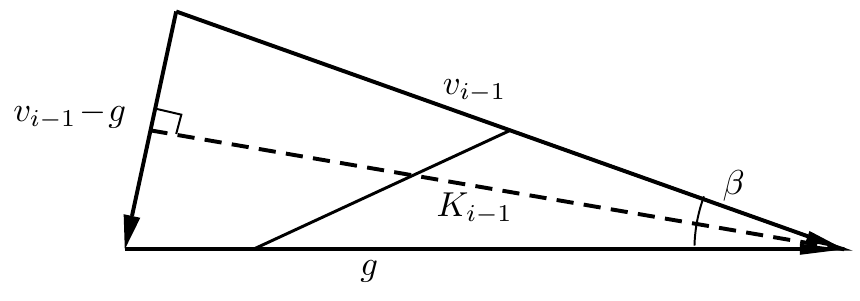}
\caption{Proof of estimate (\ref{th:H1B_error_eq10}).}
\label{fig:trianglemma1}
\end{center}
\end{figure}

It remains to estimate the term $(E)$ from (\ref{th:H1B_error_eq8}):
\begin{align}
&(E)\le\Big(\sum_{i=1}^N \frac{|\tilde K_i|}{\sin^2(\pi-\talpha_i)}\big|(\grad u_i-\grad u_{i-1})\cdot (v_{i-1}-g)\big|^2\Big)^{1/2}\label{th:H1B_error_eq9}\\
&\!+\!\!\Bigl(\sum_{i=1}^N\!\frac{|\tilde K_i|}{\sin^2(\pi\!-\!\talpha_i)} \big|(\grad U_i\!-\!\grad u_i)\cdot (v_{i-1}\!-\!g) \!-\!(\grad U_{i-1}\!-\!\grad u_{i-1})\cdot (v_{i-1}-g)\big|^2\Bigr)^{\!1/2}\!\!\!.
\nonumber
\end{align}
A simple geometric argument, cf. Figure \ref{fig:trianglemma1}, gives us, due to assumption (\ref{th:H1B_error_ass2}),
\begin{equation}
|v_{i-1}-g|=2\sin(\beta_{i-1}/2)|g|\leq 2\sin(\pi-\alpha_{i-1})\leq 2C_\alpha\sin(\pi-\talpha_{i}),
\label{th:H1B_error_eq10}
\end{equation}
since $\beta_{i-1}/2<\beta_{i-1}<\pi-\alpha_{i-1}$ and $\sin(\cdot)$ is a monotone function on $[0,\pi/2]$. Furthermore, by the definition of $u$, we have $|\grad u_i-\grad u_{i-1}|\leq Ch$. Therefore, we can estimate the first right-hand side term in (\ref{th:H1B_error_eq9}) as
\begin{equation}
\begin{split}
&\Big(\sum_{i=1}^N \frac{1}{\sin^2(\pi-\talpha_i)}\big|(\grad u_i-\grad u_{i-1})\cdot (v_{i-1}-g)\big|^2|\tilde K_i|\Big)^{1/2}\\
&\leq \Big(\sum_{i=1}^N \frac{1}{\sin^2(\pi-\talpha_i)}Ch^2\sin^2(\pi-\talpha_i)|\tilde K_i|\Big)^{1/2}\leq Ch|\B|^{1/2}\leq Ch^\alpha.
\end{split}
\label{th:H1B_error_eq11}
\end{equation}
As for the second right-hand side term in (\ref{th:H1B_error_eq9}), due to assumptions (\ref{th:H1B_error_ass2}) and (\ref{th:H1B_error_eq10}) we have
\begin{align}
&\Big(\sum_{i=1}^N\frac{|\tilde K_i|}{\sin^2(\pi\!-\!\talpha_i)} \big|(\grad U_i\!-\!\grad u_i)\!\cdot\! (v_{i-1}\!-\!g)\!-\!(\grad U_{i-1}\!-\!\grad u_{i-1})\cdot (v_{i-1}-g)\big|^2\Big)^{1/2}\nonumber\\
&\leq \Big(\sum_{i=1}^N\frac{C\sin^2(\pi-\talpha_i)}{\sin^2(\pi-\talpha_i)} \big(|\grad U_i-\grad u_i|+|\grad U_{i-1}-\grad u_{i-1}|\big)^2|\tilde K_i|\Big)^{1/2}\nonumber\\
&\leq \Big(C\sum_{i=1}^N \big|\grad U_i-\grad u_i\big|^2 C_K|K_i|\Big)^{1/2} \!\!+\Big(C\sum_{i=1}^N \big|\grad U_{i-1}-\grad u_{i-1}\big|^2 C_K|K_{i-1}|\Big)^{1/2}\nonumber\\
&\leq C\Big(\sum_{i=0}^N \big|\grad U_i-\grad u_i\big|^2 |K_i|\Big)^{1/2}  =C\Big(\sum_{i=0}^N \int_{K_i}\big|\grad U_i-\grad u_i\big|^2\dx\Big)^{1/2}\nonumber\\
&\leq C\Big(\sum_{i=0}^N \int_{K_i}\big|\grad U_i-\grad u\big|^2\dx\Big)^{1/2} +C\Big(\sum_{i=0}^N \int_{K_i}\big|\grad u-\grad u_i\big|^2\dx\Big)^{1/2}\nonumber\\
&\leq C|U-u|_1+Ch|\B|^{1/2}\leq Ch^\alpha,
\label{th:H1B_error_eq12}
\end{align}
due to (\ref{grad_ui_approx}) and (\ref{error_estimate}). 

Finally, we collect estimates (\ref{th:H1B_error_eq5})--(\ref{th:H1B_error_eq12}) of the terms $(A)$--$(E)$ and apply them in (\ref{th:H1B_error_eq4}) to estimate $\mcA_2$. Again we point out that the dominant term $\mcA_1$ appearing in the estimate of $(D)$ in (\ref{th:H1B_error_eq8}) cancels out with the corresponding term $-\mcA_1$ from (\ref{th:H1B_error_eq4}). Since all the remaining terms are of the order $O(h^\alpha)$, we have obtained the desired estimate of $\mcA_2$.
\qed

\begin{corollary}
\label{cor:H1B_error}
A necessary condition for (\ref{error_estimate}) to hold is that every band $\B$ satisfying the assumptions of Theorem \ref{th:H1B_error} must satisfy
\begin{equation}
\min_{i=1,\ldots,N}\frac{1}{\sin(\pi-\talpha_i)}\min_{i=1,\ldots,N}|\tilde K_i|^{1/2} LN^{-1/2}\leq Ch^\alpha
\label{cor:H1B_error_1}
\end{equation}
for $C$ sufficiently small.
\end{corollary}
\proof
Assume on the contrary that 
\begin{equation}
\min_{i=1,\ldots,N}\frac{1}{\sin(\pi-\talpha_i)}\min_{i=1,\ldots,N}|\tilde K_i|^{1/2} LN^{-1/2}\ge \tilde Ch^\alpha
\label{cor:H1B_error_2}
\end{equation}
for some sufficiently large $\tilde C$. The left-hand side of (\ref{cor:H1B_error_2}) is simply the estimate (up to a constant) of $\mcA_1$ from Theorem (\ref{th:H1B_error}). Then by Theorem (\ref{th:H1B_error}) we have
\begin{equation}
|u-U|_{H^1(\B)}\ge \mcA_1-\mcA_2\ge C_1\tilde Ch^\alpha-C_2h^\alpha\ge C_3h^\alpha
\label{cor:H1B_error_3}
\end{equation}
if $\tilde C$ is chosen large enough so as to dominate the term $C_2h^\alpha$ in (\ref{cor:H1B_error_3}). From (\ref{cor:H1B_error_3}) and (\ref{error_estimate}) we have
\begin{equation}
C_3h^\alpha\leq|u-U|_{H^1(\B)}\leq|u-U|_{H^1(\Omega)}\leq C(u)h^\alpha.
\nonumber
\end{equation}
For $\tilde C$ sufficiently large, $C_3$ can be made larger than $C(u)$, leading to a contradiction.
\qed

\begin{remark}
\label{rem:curvedB_2}
Here we shall comment on the case when $\Gamma$ does not lie on a straight line, cf. Remark \ref{rem:straightB}. By $g_i$ we denote the unit vector given by the edge $K_i\cap\Gamma$. Then in (\ref{lem:wdif_est}), we have $U_i^\prime=\grad U_i\cdot g_i$ and (\ref{lem:wdif_est2}) now reads
\begin{equation}
\sum_{i=2}^N |\grad U_i\cdot g_i-\grad U_{i-1}\cdot g_{i-1}|^2\ge CL^2N^{-1},
\label{lem:wdif_est2_curved}
\end{equation}
which leads to the definition of $\mcA_1$ in (\ref{th:H1B_error_eq2}) as
\begin{equation}
\mcA_1:=\Big(\sum_{i=1}^N \frac{1}{\sin^2(\pi-\talpha_i)}\big|\grad U_i\cdot g_i-\grad U_{i-1}\cdot g_{i-1}\big|^2|\tilde K_i|\Big)^{1/2},
\label{th:H1B_error_eq2_curved}
\end{equation}
for which we immediately have the estimate (\ref{th:H1B_error_eq3}) due to (\ref{lem:wdif_est2_curved}). The remaining terms in the proof of Theorem (\ref{th:H1B_error}) are the same, except for $(E)$, which is now
\begin{equation}
\Big(\sum_{i=2}^N \frac{1}{\sin^2(\pi-\talpha_i)} \big|\grad U_i\cdot (v_{i-1}-g_i)-\grad U_{i-1}\cdot (v_{i-1}-g_{i-1})\big|^2|\tilde K_i|\Big)^{1/2}.
\end{equation}
Similarly as in (\ref{th:H1B_error_eq9}),
\begin{equation}
\begin{split}
&(E)\le\Big(\sum_{i=2}^N \frac{|\tilde K_i|}{\sin^2(\pi-\talpha_i)}\big|\grad u_i\cdot (v_{i-1}-g_i)-\grad u_{i-1}\cdot (v_{i-1}-g_{i-1})\big|^2\Big)^{1/2}\\
&\!+\!\!\Big(\!\sum_{i=2}^N\!\frac{|\tilde K_i|}{\sin^2(\pi\!-\!\talpha_i)} \big|(\grad U_i\!-\!\grad u_i)\!\cdot\! (v_{i-1}\!-\!g_i)\!-\!(\grad U_{i-1}\!-\!\grad u_{i-1})\!\cdot\! (v_{i-1}\!-\!g_{i-1})\big|^2\Big)^{\!1/2}
\end{split}
\label{main_est:eq8_curved}
\end{equation}
where the second term can be estimated as (\ref{th:H1B_error_eq12}), since $|v_{i}-g_{i}|\leq 2C_\alpha\sin(\pi-\talpha_{i})$ and similarly for $v_{i-1}-g_{i-1}$ by Figure \ref{fig:trianglemma1}. The first term in (\ref{main_est:eq8_curved}) can be rewritten using the triangle inequality 
\begin{equation}
\begin{split}
&\big|\grad u_i\cdot (v_{i-1}-g_i)-\grad u_{i-1}\cdot (v_{i-1}-g_{i-1})\big|\\
&\leq \big|(\grad u_i-\grad u_{i-1})\cdot (v_{i-1}-g_{i-1})\big| +\big|\grad u_{i}\cdot (g_{i-1}-g_i)\big|\\
&\leq Ch\sin(\pi-\talpha_{i})+C|g_{i-1}-g_i|.
\end{split}
\end{equation}
In order to eliminate the denominator $\sin(\pi-\talpha_{i})$ in (\ref{main_est:eq8_curved}) and to obtain an $O(h^{\alpha})$ estimate for $(E)$, one must e.g. assume that $|g_{i-1}-g_i|\leq C\sin(\pi-\talpha_{i}) h^{\alpha-1/2}$ along with $|\B|\leq Ch$. Another possibility is to estimate $|\B|\leq CLh\max_{i=1,\ldots,N} \sin(\pi-\talpha_{i})$, which (under additional assumptions similar to (\ref{th:H1B_error_ass2}) eliminates the unpleasant denominator) and assume $|g_{i-1}-g_i|\leq Ch^{\alpha-1/2}$. This is possible for one band. For multiple bands as in Section \ref{sec:main_est:mult_band} one must assume e.g. $|g_{i-1}-g_i|\leq C\sin(\pi-\talpha_{i}) h^{\alpha}$. We note that without any assumptions it naturally holds $|g_{i-1}-g_i|\leq \sin(\pi-\talpha_{i})+\sin(\pi-\alpha_{i-1})+\sin(\pi-\alpha_{i})$, hence we assume an additional factor of e.g. $O(h^{\alpha-1/2})$. Overall, these variants correspond to $\B$ not lying on a straight line but with a restriction on how much it `bends'.
\end{remark}

\subsection{Examples and counterexamples -- one band}
\label{subsec:necessary_examples}
\begin{figure}[t]
\begin{center}
\includegraphics[scale=0.7,clip]{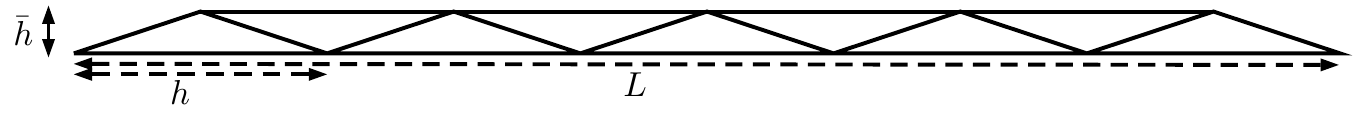}
\caption{A band consisting of identical isosceles triangles.}
\label{fig:Band2}
\end{center}
\end{figure}

Here we present several specific applications of Corollary \ref{cor:H1B_error}. For simplicity, we assume that the band $\B\subset\Th$ consists of identical isosceles triangles with base length $h$ and height $\bar h$, cf. Figure \ref{fig:Band2}. Such a regular structure is by no means required by Theorem (\ref{th:H1B_error}), however the evaluation of all quantities in estimate (\ref{th:H1B_error_est2}) is very simple in this case and also the specific examples can be directly compared to the results of \cite{Babuska-Aziz}, \cite{Oswald}.

We assume that $\B=(\cup_{i=0}^N K_i)\cup(\cup_{i=1}^N\tilde K_i)$ has length $L=(N+1)h$ for some $N\in\mathbb{N}$. Then
\begin{equation}
N^{-1/2}\ge(N+1)^{-1/2}=h^{1/2}L^{-1/2}.
\label{nec_ex1}
\end{equation}
For all $i$, we have $|\tilde K_i|=\tfrac{1}{2}h\bar{h}$. Finally, we need to estimate $\sin(\pi-\tilde\alpha_i)$. Since all $\tilde K_i$ are identical, we denote simply $\talpha_i=\talpha$.

\begin{lemma}
\label{lem:sin}
We have
\begin{equation}
\frac{1}{\sin(\pi-\talpha)}>\frac{h}{4\bar h}.
\label{lem:sin_est}
\end{equation}
\end{lemma}
\proof
\begin{figure}[t]
\begin{center}
\includegraphics[scale=0.7,clip]{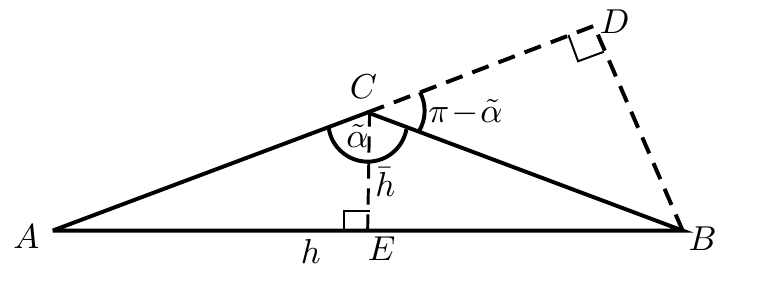}
\caption{Proof of Lemma \ref{lem:sin}.}
\label{fig:trianglemma2}
\end{center}
\end{figure}
Let $\tilde{K}$ be the triangle $ABC$ in Figure \ref{fig:trianglemma2}. Triangles $ABD$ and $ACE$ are similar, therefore 
\begin{equation}
\frac{|CB|\sin(\pi-\talpha)}{\bar h}=\frac{|DB|}{\bar{h}}=\frac{h}{|AC|}<2,
\label{lem:sin2}
\nonumber
\end{equation}
since $\tfrac{1}{2}h<|AC|$. Hence $|CB|\sin(\pi-\talpha)<2\bar h$. Since $\tfrac{1}{2}h<|CB|$, we get $\tfrac{1}{2}h\sin(\pi-\talpha)<2\bar h$, which is (\ref{lem:sin_est}).
\qed

\medskip
Using (\ref{nec_ex1}) and (\ref{lem:sin_est}), the left-hand side of the necessary condition (\ref{cor:H1B_error_1}) from Corollary \ref{cor:H1B_error} can be estimated as
\begin{equation}
\begin{split}
\min_{i=1,\ldots,N}&\frac{1}{\sin(\pi-\talpha_i)}\min_{i=1,\ldots,N}|\tilde K_i|^{1/2} LN^{-1/2}\\ &\ge \frac{h}{4\bar h}(\tfrac{1}{2}h\bar h)^{1/2}L(h^{1/2}L^{-1/2})=\frac{1}{4\sqrt{2}}h^2\bar h^{-1/2}L^{1/2}.
\end{split}
\label{nec_ex2}
\end{equation}
Corollary \ref{cor:H1B_error} along with (\ref{nec_ex2}) then gives us the necessary condition for $O(h^\alpha)$ convergence in the form
\begin{equation}
\frac{1}{4\sqrt{2}}h^2\bar h^{-1/2}L^{1/2}\leq Ch^\alpha
\label{nec_ex3}
\end{equation}
for $C$ sufficiently small. Expressing $\bar h$, we get the necessary condition
\begin{equation}
\bar h\ge Ch^{4-2\alpha}L
\label{nec_ex3a}
\end{equation}
for $C$ sufficiently small.

\subsubsection{Bands of minimal length}
\label{subsubsec:bands_min_len}
Theorem \ref{th:H1B_error} requires that $L\ge C_Lh^{2\alpha/5}$ for $C_L$ sufficiently large. If we take $L=C_Lh^{2\alpha/5}$ in (\ref{nec_ex3a}), i.e. the minimal possible length of $\B$, we get the following.
\begin{lemma}
\label{lem:nec1}
Let $\B\subset\Th$ be a band of length $L=C_Lh^{2\alpha/5}$ with the regular geometry considered above, then a necessary condition for (\ref{error_estimate}) to hold, i.e. for $O(h^\alpha)$ convergence is
\begin{equation}
\bar h\ge Ch^{4-8\alpha/5}
\label{nec_ex4}
\end{equation}
for $C$ sufficiently large.
\end{lemma}

Of course, $L=C_Lh^{2\alpha/5}$ need not be satisfied exactly for all $h$, we can have $L\sim C_Lh^{2\alpha/5}$.  As a straightforward consequence of Lemma (\ref{lem:nec1}), we get the following counterexamples:

\medskip
\noindent{\bf Counterexample 1:} Let $\{\Th\}_{h\in(0,h_0)}$ be such that a sequence of bands $\B_{k}\subset\mathcal{T}_{h_k}, k\in\mathbb{N}\,$, exists for $h_k\to 0$ as $k\to\infty$, with the regular geometry considered above. Let the length and shape parameters of each $\B\in\{\B_{k}\}_{k\in\mathbb{N}}$ satisfy
\begin{equation}
L= C_Lh^{2\alpha/5},\quad \bar h=o( h^{4-8\alpha/5}),
\label{nec:counter:1}
\end{equation}
for $C_L$ large enough, then (\ref{error_estimate}) does not hold for $u$ given by (\ref{u_def_quadfunc}) for any $C(u)\ge 0$, i.e. the FEM cannot have $O(h^\alpha)$ convergence on $\{\Th\}_{h\in(0,h_0)}$.

We note that since we essentially consider only the error on $\B$, the counterexample is independent of the properties of the rest of the triangulation, i.e. $\Th\setminus\B$ can be as `nice' as possible, structured, uniform, etc. 

\medskip
As special cases of Counterexample 1, we take the most interesting values $\alpha=1$ and $\alpha=0$.

\medskip
\noindent{\bf Counterexample 2:} As in Counterexample 1, let there be an infinite sequence of bands $\B\subset\Th$ with $h\to 0$ (we omit the subscript $k$ for simplicity) such that
\begin{equation}
L= C_Lh^{2/5},\quad \bar h=o( h^{12/5})=o(h^{2.4}),
\label{nec:counter:2}
\end{equation}
then the FEM cannot have $O(h)$ convergence on $\{\Th\}_{h\in(0,h_0)}$ for general $u$.

\medskip
\noindent{\bf Counterexample 3:} As in Counterexample 1, let there be an infinite sequence of bands $\B\subset\Th$ with $h\to 0$ such that
\begin{equation}
L\ge C_L,\quad \bar h=o(h^4),
\label{nec:counter:3}
\end{equation}
for any fixed $C_L>0$. Then the FEM does not converge in $H^1(\Omega)$ on $\{\Th\}_{h\in(0,h_0)}$ for general $u$. 

\proof
Assume on the contrary that $|u-U|_1\to 0$ as $h\to 0$. Choose $\varepsilon>0$ arbitrary but fixed. Then for all $h$ sufficiently small $|u-U|_1\leq\varepsilon$, i.e estimate (\ref{error_estimate}) holds with $\alpha=0, C(u)=\varepsilon$. From Theorem \ref{th:H1B_error}, if $L\ge C_L$ sufficiently large w.r.t. $C(u)=\varepsilon$, we have
\begin{equation}
|u-U|_1\ge \mathcal{A}_1-\mathcal{A}_2\to\infty \quad\text{as } h\to 0,
\label{nec:counter:3a}
\end{equation}
since $|\mathcal{A}_2|\leq C_2$, a constant, and due to (\ref{nec_ex2}) and the assumption $\bar{h}=o(h^4)$, $\mathcal{A}_1\to\infty$ as $h\to 0$. Estimate (\ref{nec:counter:3a}) is in contradiction with the assumption $|u-U|_1\leq\varepsilon$ for all $h$ sufficiently small. Hence $|u-U|_1\not\to 0$.

We note that to be able to use Theorem \ref{th:H1B_error}, $C_L$ must be sufficiently large with respect to the fixed $\varepsilon$. However, for the sake of the proof by contradiction, $\varepsilon$ can be taken arbitrarily small. By Remark \ref{rem:C_suff_large}, for a given lower bound $C_L>0$ for $L$ in (\ref{nec:counter:3}), we can always choose $\varepsilon$ small enough so that this particular $C_L$ satisfies the assumptions of Theorem \ref{th:H1B_error}. The ability to take $C_L$ arbitrarily small but fixed is important, since for e.g. $\varepsilon=1$ a band of length $C_L$ might not even fit into $\Omega$.
\qed


\subsubsection{Bands of constant length}
\label{subsubsec:bands_cons_len}
Taking $L\ge C_L$ in Counterexample 3 was necessary for $\alpha=0$ due to the condition $L\ge C_Lh^{2\alpha/5}$ in Theorem \ref{th:H1B_error}. On the other hand, one can take $L\sim C_L$ for all $\alpha\in[0,1]$ instead of the minimal $L$ from Section \ref{subsubsec:bands_min_len}. Taking $L=C_L$ in (\ref{nec_ex3a}), we get:

\begin{lemma}
\label{lem:nec2}
Let $\B\subset\Th$ be a band of length $L=C_L$ with the regular geometry considered above, then a necessary condition for $O(h^\alpha)$ convergence is
\begin{equation}
\bar h\ge Ch^{4-2\alpha}
\label{nec_ex4a}
\end{equation}
for $C$ sufficiently large.
\end{lemma}

\medskip
\noindent{\bf Counterexample 4:}
Let there be an infinite sequence of bands $\B\subset\Th$ with $h\to 0$ such that
\begin{equation}
L= C_L,\quad \bar h=o( h^{4-2\alpha}),
\label{nec:counter:4}
\end{equation}
then the FEM cannot have $O(h^\alpha)$ convergence on $\{\Th\}_{h\in(0,h_0)}$ for general $u$.

\medskip
Again, as a special case, we take the most interesting value $\alpha=1$. The case $\alpha=0$ is covered by Counterexample 3.

\medskip
\noindent{\bf Counterexample 5:}
Let there be an infinite sequence of bands $\B\subset\Th$ with $h\to 0$ such that
\begin{equation}
L= C_L,\quad \bar h=o( h^{2}),
\label{nec:counter:5}
\end{equation}
then the FEM cannot have $O(h)$ convergence on $\{\Th\}_{h\in(0,h_0)}$ for general $u$.

\subsection{Main estimate -- multiple bands}
\label{sec:main_est:mult_band}
Up to now, we dealt with estimates on one band only. However, more can be gained by considering the case when there are multiple bands $\B$ in $\Th$. Specifically, we shall consider $N_\B$ disjoint bands $\{\B_b\}_{b=1}^{N_\B}\subset\Th$. We use the notation $\BB=\cup_{b=1}^{N_\B}\B$. Of course, Lemmas \ref{lem:angles}--\ref{lem:wdif} remain valid for each  $\B_b$. It is the analysis of Section \ref{sec:main_est:one_band} that needs to be modified to combine these individual estimates together. For example, simply summing the necessary condition (\ref{cor:H1B_error_2}) over all bands would not yield anything new. However, estimating the error on the whole union $\cup_{b=1}^{N_\B}\B_b$ gives the following stronger result of Theorem \ref{th:mb:H1B_error}.

Since we now consider multiple bands, we need to take this into account in the notation, which is analogous to that of Section \ref{sec:main_est:one_band}. Let $\cup_{b=1}^{N_\B}\B_b\subset\Th$. For $b=1,\ldots,N_\B$, let $\B_b=(\cup_{i=0}^{N_b} K_i^b)\cup(\cup_{i=1}^{N_b} \tilde K_i^b)$ be the band with edge $\Gamma_b$. We denote $U_i^b=U|_{K_i^b}$ and $\tilde U_i^b=U|_{\tilde K_i^b}$.
By $g_b$ we denote the unit vector in the direction of $\Gamma^b$ and by $L_b$ the length of $\B_b$. 

Similarly as in (\ref{lem:wdif_est2}), we can reformulate Lemma \ref{lem:wdif} for $\B_b$ as
\begin{equation}
\sum_{i=1}^{N_b} \big|(\grad U_i^b-\grad U_{i-1}^b)\cdot g_b\big|^2\ge CL_b^2N_b^{-1}.
\label{lem:mbwdif_est2}
\end{equation}
 
Again, we denote $\grad u_i^b=\grad u(\mcC_i^b)$, where $\mcC_i^b$ is the centroid of $K_i^b$ and similarly for $\tilde K^b_i$. We have approximation properties similar to (\ref{grad_ui_approx}). By $\alpha_i^b,\talpha_i^b$ we denote the maximum angles of $K_i^b$ and $\tK_i^b$, respectively.

\begin{theorem}
\label{th:mb:H1B_error}
Let $u$ be given by (\ref{u_def_quadfunc}) and let $U\in X_h$ satisfy the error estimate (\ref{error_estimate}). Let the set of bands $\BB=\cup_{b=1}^{N_\B}\B_b\subset\Th$ with $\B_b=(\cup_{i=0}^{N_b} K_i^b)\cup(\cup_{i=1}^{N_b} \tilde K_i^b)$ be such that there exist constants $C_K, C_\alpha$ independent of $h$ such that for all $b=1,\ldots,N_\B$ and $i=1,\ldots,N_b$
\begin{equation}
\begin{split}
|\tilde K_i^b|&\leq C_K|K_i^b|,\\
|\tilde K_i^b|&\leq C_K|K_{i-1}^b|,\\
\sin(\pi-\alpha_{i-1}^b)&\leq C_\alpha\sin(\pi-\tilde\alpha_{i}^b).
\label{th:mb:H1B_error_ass2}
\end{split}
\end{equation}
Finally, let $L_b\ge C_Lh^{2\alpha/5}$, for all $b=1,\ldots,N_\B$, where $C_L$ is sufficiently large. Then 
\begin{equation}
|u-U|_{H^1(\BB)}\ge \mcA_1-\mcA_2,
\label{th:mb:H1B_error_est1}
\end{equation}
where 
\begin{equation}
\begin{split}
\mcA_1&\ge C_1\Big(\sum_{b=1}^{N_\B}\min_{i=1,\ldots,N_b} \frac{1}{\sin^2(\pi-\talpha_i^b)}\min_{i=1,\ldots,N_b}|\tilde K_i^b| L_b^2N^{-1}_b\Big)^{1/2},\\
|\mcA_2|&\leq C_2h^\alpha.
\end{split}
\label{th:mb:H1B_error_est2}
\end{equation}
\end{theorem}
\proof
The proof is completely analogous to that of Theorem \ref{th:H1B_error}, therefore we only point out the main differences. Similarly as in (\ref{th:H1B_error_eq1}), we consider the error only on the elements $\tilde K_i^b$:
\begin{equation}
\begin{split}
|u&-U|_{H^1(\BB)}\ge \Big(\sum_{b=1}^{N_\B}\sum_{i=1}^{N_b}\int_{\tilde K_i^b}|\grad u-\grad U|^2\dx\Big)^{1/2} =\mcA_1-\mcA_2,
\end{split}
\label{th:mb:H1B_error_eq1}
\end{equation}
where the splitting of the error is given by
\begin{equation}
\begin{split}
\mcA_1&:=\Big(\sum_{b=1}^{N_\B}\sum_{i=1}^{N_b} \frac{1}{\sin^2(\pi-\talpha_i^b)}\big|(\grad U_i^b-\grad U_{i-1}^b)\cdot g_b\big|^2|\tilde K_i^b|\Big)^{1/2},\\
\mcA_2&:=\Big(\sum_{b=1}^{N_\B}\sum_{i=1}^{N_b}\int_{\tilde K_i^b}|\grad u-\grad U|^2\dx\Big)^{1/2}-\mcA_1.
\label{th:mb:H1B_error_eq2}
\end{split}
\end{equation}
From (\ref{lem:wdif_est2}) we immediately get the first estimate of (\ref{th:mb:H1B_error_est2}). 

Similarly as in (\ref{th:H1B_error_eq4}), we estimate
\begin{align}
\mcA_2&\leq \Big(\sum_{b=1}^{N_\B}\sum_{i=1}^{N_b}\int_{\tilde K_i^b}|\grad u-\grad \tilde u_i^b|^2\dx\Big)^{1/2}+
\Big(\sum_{b=1}^{N_\B}\sum_{i=1}^{N_b}|\grad \tilde u_i^b-\grad u_i^b|^2|\tilde K_i^b|\Big)^{1/2}\nonumber\\
&+\!\Big(\!\sum_{b=1}^{N_\B}\sum_{i=1}^{N_b}|\grad u_i^b-\grad U_i^b|^2|\tilde K_i^b|\Big)^{\!1/2} \!\!+\!\Big(\!\sum_{b=1}^{N_\B}\sum_{i=1}^{N_b}|\grad U_i^b-\grad \tilde U_i^b|^2|\tilde K_i^b|\Big)^{\!1/2}\!\!\!-\!\mcA_1\nonumber\\
&=:(A)+(B)+(C)+(D)-\mcA_1
\label{th:mb:H1B_error_eq4}
\end{align}
To demonstrate the differences between the proof of Theorems \ref{th:H1B_error} and \ref{th:mb:H1B_error}, we estimate only e.g. term $(C)$. Similarly as in (\ref{th:H1B_error_eq7}), we have
\begin{align}
(C)&\leq \Big(C_K\!\sum_{b=1}^{N_\B}\sum_{i=1}^{N_b}\int_{K_i^b}\!|\grad u_i^b-\grad u|^2\dx\Big)^{\!1/2} \!\!+\!\Big(C_K\!\sum_{b=1}^{N_\B}\sum_{i=1}^{N_b}\int_{K_i^b}\!|\grad u-\grad U|^2\dx\Big)^{\!1/2}\nonumber\\
&\leq Ch|\BB|^{1/2}+C_K|u-U|_{H^1(\BB)}\leq Ch^\alpha
\label{th:mb:H1B_error_eq7}
\end{align}
due to (\ref{grad_ui_approx}) and (\ref{error_estimate}). This is the main difference from the proof of Theorem \ref{th:H1B_error}: the estimation $|\BB|^{1/2}\leq|\Omega|^{1/2}$ instead of $|\B|^{1/2}\leq|\Omega|^{1/2}$ and using the error estimate (\ref{error_estimate}) on the whole union $\BB$ instead of individual bands $\B_b$ and then summing over $b$, which would yield the undesired estimate $N_\B Ch^\alpha$. All the other terms can be estimated as in Theorem \ref{th:H1B_error} with this in mind.
\qed

\medskip
From Theorem \ref{th:mb:H1B_error}, we get a necessary condition for $O(h^\alpha)$ convergence similar to Corollary \ref{cor:H1B_error}.

\begin{corollary}
\label{cor:mb:H1B_error}
A necessary condition for (\ref{error_estimate}) to hold is that every set of bands $\BB\subset\Th$ satisfying the assumptions of Theorem \ref{th:H1B_error} must satisfy
\begin{equation}
\Big(\sum_{b=1}^{N_\B}\min_{i=1,\ldots,N_b} \frac{1}{\sin^2(\pi-\talpha_i^b)}\min_{i=1,\ldots,N_b}|\tilde K_i^b| L_b^2N^{-1}_b\Big)^{1/2}\leq Ch^\alpha
\label{cor:mb:H1B_error_1}
\end{equation}
for $C$ sufficiently small.
\end{corollary}

\subsection{Examples and counterexamples -- multiple bands}
\label{subsec:necessary_examples:mb}
\begin{figure}[t]
\begin{center}
\includegraphics[scale=0.4,clip]{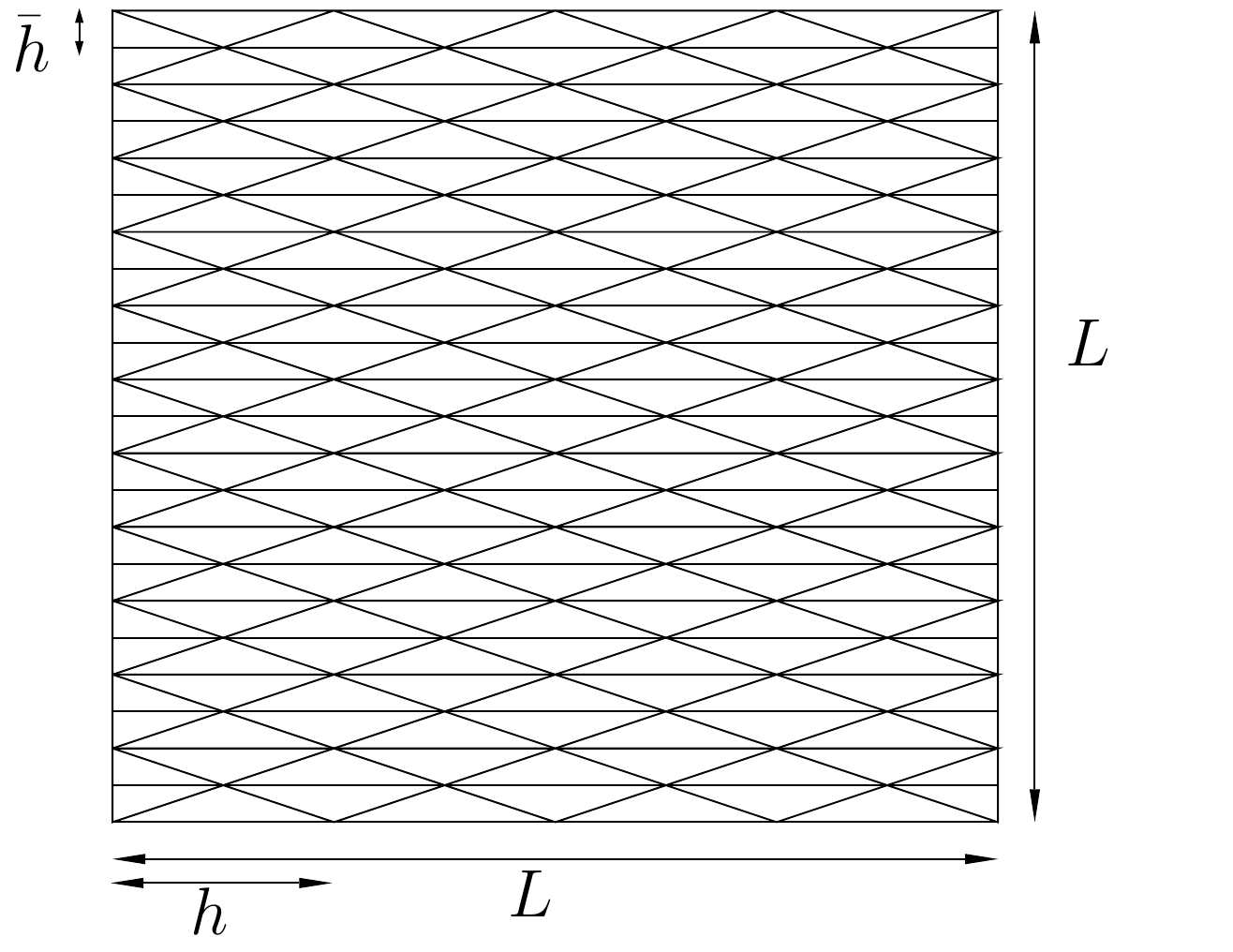}
\caption{A Babu\v{s}ka-Aziz-type set of bands $\BB$ consisting of identical isosceles triangles}
\label{fig:Babuska-Aziz}
\end{center}
\end{figure}

Here we present some applications of Corollary \ref{cor:mb:H1B_error}. As in Section \ref{subsec:necessary_examples}, we assume for simplicity that the bands $\B_b\subset\Th$ consists of identical isosceles triangles with base length $h$ and height $\bar h$. Moreover, we assume that the union of clusters $\BB$ has the structure of the Babu\v{s}ka-Aziz counterexample, cf. Figure \ref{fig:Babuska-Aziz}. Again we note that such a regular structure of the bands is not required by Theorem (\ref{th:H1B_error}) and the bands do not need to lie next to each other as in Figure \ref{fig:Babuska-Aziz}. We consider such a simple geometry to make the estimates as clear as possible and directly comparable to the optimal results of \cite{Oswald} on the Babu\v{s}ka-Aziz counterexample.

Similarly as in Section \ref{subsec:necessary_examples}, we have the estimates
\begin{equation}
N_b^{-1}\ge hL^{-1},\quad |\tilde K_i^b|=\tfrac{1}{2}h\bar{h}, \quad \frac{1}{\sin^2(\pi-\talpha)}>\frac{h^2}{16\bar h^2}.
\label{nec_ex0:mb}
\end{equation}
Finally, due to the structure of $\BB$ considered, we have $N_\B=L\bar h^{-1}$.

Then the left-hand side of the necessary condition (\ref{cor:mb:H1B_error_1}) from Corollary \ref{cor:mb:H1B_error} can be estimated as
\begin{equation}
\begin{split}
\Big(&\sum_{b=1}^{N_\B}\min_{i=1,\ldots,N_b} \frac{1}{\sin^2(\pi-\talpha_i^b)}\min_{i=1,\ldots,N_b}|\tilde K_i^b| L_b^2N^{-1}_b\Big)^{1/2}\\
&\ge \Big(N_\B\frac{h^2}{16\bar h^2}(\tfrac{1}{2}h\bar h)L^2(hL^{-1})\Big)^{1/2} =\Big((L\bar h^{-1})\frac{1}{32}h^4\bar h^{-1}L\Big)^{1/2}=\frac{1}{4\sqrt{2}}h^2\bar h^{-1}L.
\label{nec_ex2:mb}
\end{split}
\end{equation}
Corollary \ref{cor:mb:H1B_error} along with (\ref{nec_ex2:mb}) then gives us the necessary condition for $O(h^\alpha)$ convergence in the form
\begin{equation}
\frac{1}{4\sqrt{2}}h^2\bar h^{-1}L\leq Ch^\alpha
\label{nec_ex3:mb}
\end{equation}
for $C$ sufficiently small. Expressing $\bar h$ as in (\ref{nec_ex3a}), we get the necessary condition
\begin{equation}
\bar h\ge Ch^{2-\alpha}L
\label{nec_ex3:mba}
\end{equation}
for $C$ sufficiently small.

\subsubsection{Multiple bands of minimal length}
\label{subsubsec:mb:bands_min_len}
As in Lemma \ref{lem:nec1}, we can take the minimal length $L=C_Lh^{2\alpha/5}$ in (\ref{nec_ex3:mba}) to obtain:

\begin{lemma}
\label{lem:mb:nec1}
Let $\BB\subset\Th$ be a set of bands of length $L=C_Lh^{2\alpha/5}$ with the regular geometry considered above, then a necessary condition for $O(h^\alpha)$ convergence is
\begin{equation}
\bar h\ge Ch^{2-3\alpha/5}
\label{nec_ex4:mb}
\end{equation}
for $C$ sufficiently large.
\end{lemma}

\medskip
\noindent{\bf Counterexample 6:} Let there be an infinite sequence of sets of bands $\BB\subset\Th$ with $h\to 0$ (again omitting the subscript $k$ for simplicity) such that
\begin{equation}
L= C_Lh^{2/5},\quad \bar h=o( h^{7/5})=o(h^{1.4}),
\label{nec:counter:6}
\end{equation}
then the FEM cannot have $O(h)$ convergence on $\{\Th\}_{h\in(0,h_0)}$ for general $u$.

\medskip
\noindent{\bf Counterexample 7:} Let there be an infinite sequence of sets of bands $\BB\subset\Th$ with $h\to 0$ such that
\begin{equation}
L\ge C_L,\quad \bar h=o(h^2),
\label{nec:counter:7}
\end{equation}
for any fixed $C_L>0$. Then the FEM does not converge in $H^1(\Omega)$ on $\{\Th\}_{h\in(0,h_0)}$ for general $u$.

\subsubsection{Multiple bands of constant length}
\label{subsubsec:mb:bands_cons_len}
By taking  $L\ge C_L$ in (\ref{nec_ex3:mba}), we get:

\begin{lemma}
\label{lem:mb:nec2}
Let $\BB\subset\Th$ be a set of bands of length $L=C_L$ with the regular geometry considered above, then a necessary condition for $O(h^\alpha)$ convergence is
\begin{equation}
\bar h\ge Ch^{2-\alpha}
\label{nec_ex5:mb}
\end{equation}
for $C$ sufficiently large.
\end{lemma}

\medskip
\noindent{\bf Counterexample 8:} Let there be an infinite sequence of sets of bands $\BB\subset\Th$ with $h\to 0$ such that
\begin{equation}
L= C_L,\quad \bar h=o(h),
\label{nec:counter:8}
\end{equation}
then the FEM cannot have $O(h)$ convergence on $\{\Th\}_{h\in(0,h_0)}$ for general $u$.

\begin{remark}
\label{rem:oswald}
We note that the result of Lemma \ref{lem:mb:nec2} is optimal when applied to the Babu\v ska-Aziz counterexample, since in \cite{Oswald} it is proven that in this case (using our notation)
\begin{equation}
|u-U|_1\sim \min(1,h^2/\bar h).
\label{rem:oswald1}
\end{equation}
If we want to obtain a condition on $\bar h$ for $O(h^\alpha)$-convergence from (\ref{rem:oswald1}), we get exactly (\ref{nec_ex5:mb}).
\end{remark}

In the Babu\v ska-Aziz counterexample, \cite{Babuska-Aziz}, the unit square is divided into a triangulation consisting of bands of the considered type, as in Figure \ref{fig:Babuska-Aziz}. If we consider only one band in $\Th$, as in Section \ref{subsec:necessary_examples}, we need stronger assumptions to produce the corresponding counterexamples. For example, comparing Counterexamples 5 and 8, if there is only one band we need $\bar h=o(h^2)$ for counterexamples to $O(h)$ convergence with $L=C_L$, however only $\bar h=o(h)$ is needed if we have multiple bands.

\begin{remark}
\label{rem:nec:angles}
Another possibility how to view the necessary conditions (\ref{nec_ex3}) and (\ref{nec_ex3:mb}) is by expressing these conditions not for $\bar h$, but for $\pi-\tilde\alpha_i$. Trivially, $\pi-\tilde\alpha_i>\sin(\pi-\tilde\alpha_i)$. Moreover, similarly as in Lemma \ref{lem:sin} it is possible to prove that $\sin(\pi-\tilde\alpha_i)>\bar h/h$. Combining these results with (\ref{nec_ex3}) gives us the necessary condition for $O(h^\alpha)$-convergence
\begin{equation}
\pi-\tilde\alpha_i\ge Ch^{3-2\alpha}L
\nonumber
\end{equation}
for all $i=1,\ldots,N_b$ and for $C$ sufficiently small. Similarly, from (\ref{nec_ex3:mba}), we get the necessary condition 
\begin{equation}
\pi-\tilde\alpha_i^b\ge Ch^{1-\alpha}L
\nonumber
\end{equation}
for all $i=1,\ldots,N_b$, $b=1,\ldots,N_\B$ and for $C$ sufficiently small. These angle conditions are mentioned in the Introduction.
\end{remark}

\section{A sufficient condition for FEM convergence}
\label{sec:sufficient}
A basic tool in FEM error analysis is C\' ea's lemma, cf. \cite{Ciarlet}. In fact, this seems to be the only general tool in FEM analysis that is completely independent of the geometry of $\Th$.
\begin{lemma}[C\' ea]
\begin{equation}
|u-U|_1= \inf_{v_h\in V_h}|u-v_h|_1.
\label{Cea}
\end{equation}
\end{lemma}
In a priori error analysis, one typically proceeds by taking $v_h$ in (\ref{Cea}) to be the Lagrange interpolation of $u$. This choice leads to the well known maximum angle condition for $O(h)$ convergence and the newer, lesser known circumradius condition \cite{Kobayashi}, \cite{Rand} for convergence of the FEM. We note that by using C\' ea's lemma, we are again dealing only with the approximation properties of the space $V_h$ (or $X_h$), the specific form of the FEM scheme is not taken into account, cf. Remark \ref{rem:approx}.

We proceed as follows. Since Lagrange interpolation gives suitable estimates only on elements satisfying the maximum angle condition, we will use a modified version of Lagrange interpolation on the remaining elements. This is done by dropping the interpolation condition at the maximum angle vertex of $K$. The question is then how to connect such interpolants together continuously. For this purpose we introduce the concept of a correction function which is used to construct a globally continuous function from $V_h$. Basically, we modify the interpolated function $u$ locally so that this new function $\tilde u$ is $O(h)$-close to $u$ and can be interpolated by the standard Lagrange procedure with an $O(h)$ error on $\Th$.

\subsection{Standard Lagrange interpolation}
Let $K\in \Th$ have vertices $A_K,B_K,C_K$, where $A_K$ will denote the maximum angle vertex of $K$ throughout this section, cf. Figure \ref{fig:triang}. For each element $K\in\Th$, we seek $\vhK\in P^1(K)$ such that the following conditions hold:
\begin{equation}
\begin{split}
\label{eq:th1Kdef}
\vhK(A_K)&=u(A_K),\\
\vhK(B_K)&=u(B_K),\\
\vhK(C_K)&=u(C_K),
\end{split}
\end{equation}
For the $H^1$-error of such an interpolation, we have the following result, \cite{Babuska-Aziz}, \cite{Krizek}:
\begin{lemma}
\label{lem:th1est}
Let $\alpha_0<\pi$ be fixed. Let $K\in\Th$ with maximum angle $\alpha_K\leq\alpha_0$. Then for $\vhK$ defined by (\ref{eq:th1Kdef}), we have
\begin{equation}
\big|u-\vhK\big|_{H^1(K)}\leq  C(\alpha_0)h_K|u|_{H^2(K)},
\label{lem:th1est:1}
\end{equation}
where $C(\alpha_0)$ is independent of $u, K$.
\end{lemma}

We now define the (global) Lagrange interpolation of $u$ as the function $\Pi_h u\in V_h$ such that $(\Pi_h u)|_K=\vhK$ for all $K\in\Th$. By taking $\Pi_h u$ in C\' ea's lemma, we get the \emph{maximum angle condition}.
\begin{theorem}[Maximum angle condition, \cite{Babuska-Aziz}, \cite{Krizek}]
\label{lem:max_angle}
Let $\alpha_K\leq\alpha_0<\pi$ for all $K\in\Th, h\in(0,h_0)$. Then 
\begin{equation}
|u-U|_1\leq  C(\alpha_0)h|u|_2.
\nonumber
\end{equation}
\end{theorem}

The maximum angle condition is a condition for $O(h)$ convergence, which is the subject of Section \ref{sec:sufficient}. The general case of $O(h^\alpha)$ convergence (as in Section \ref{sec:necessary}) is a simple extension treated in Remark \ref{rem:O_h_alpha}. In this case, the maximum angle condition is too strong an assumption. Recently, a generalization of Lemma \ref{lem:th1est} and the maximum angle condition was derived, cf. \cite{Kobayashi}, \cite{Rand}.

\begin{lemma}
\label{lem:circumradius}
Let $R_K\leq 1$ be the circumradius of the triangle $K\in\Th$. Then
\begin{equation}
|u-\vhK|_{H^1(K)}\leq  CR_K|u|_{H^2(K)}.
\label{lem:circumradius_est}
\end{equation}
\end{lemma}

\begin{theorem}[Circumradius estimate, \cite{Kobayashi}, \cite{Rand}]
\label{thm:circumradius}
\begin{equation}
|u-U|_1\leq  C\max_{K\in\Th}R_K|u|_2.
\nonumber
\end{equation}
\end{theorem}

\noindent The \emph{circumradius condition} of \cite{Kobayashi} then reads
\begin{equation}
\max_{K\in\Th}R_K\to 0
\nonumber
\end{equation}
and is a condition on convergence (not $O(h)$ convergence) of the FEM by Theorem \ref{thm:circumradius}. We note that the law of sines states
\begin{equation}
\frac{h_K}{\sin{\alpha_K}}=2R_K.
\nonumber
\end{equation}
If we substitute this expression into (\ref{lem:circumradius_est}), we get $O(h)$ convergence if and only if the denominator $\sin\alpha_K$ is uniformly bounded away from zero for all $K$, which is exactly the maximum angle condition. Therefore, as far as $O(h)$ convergence is concerned, Lemmas \ref{lem:th1est} and \ref{lem:circumradius} are equivalent.

In \cite{Hannukainen}, it is shown that the maximum angle condition is not necessary for $O(h)$ convergence of the finite element method, i.e. $\Th$ can contain elements whose maximum angle degenerates to $\pi$ with respect to $h\to 0$. The argument in \cite{Hannukainen} can also be used to show that the circumradius condition is not necessary for convergence, i.e $\Th$ can contain arbitrarily `bad' elements. In the following section, we modify the Lagrange interpolation on such elements to obtain a generalization of the maximum angle and circumradius conditions.

\subsection{Modified Lagrange interpolation}
\label{subsec:modif_lagr}
\begin{figure}[t]
\centering
\includegraphics[height=1.9cm]{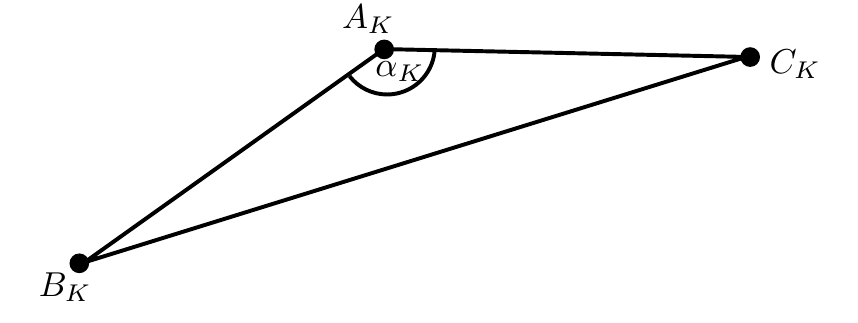}
\includegraphics[height=2.1cm]{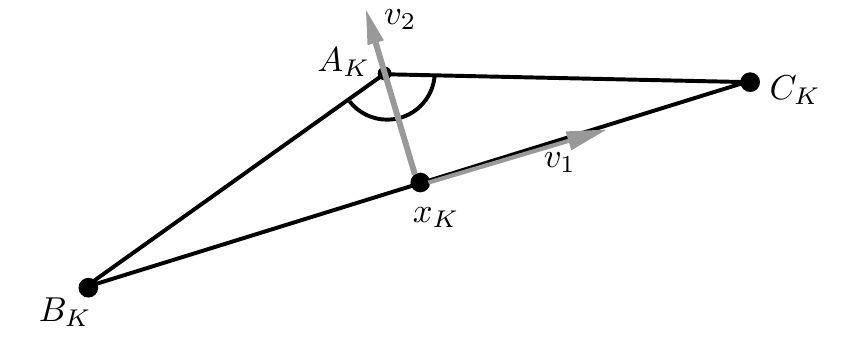}
\caption{Lagrange interpolation conditions -- left. Modified Lagrange interpolation -- right.}
\label{fig:triang}
\end{figure}

We split $\Th$ into two subsets of elements: those satisfying the maximum angle condition and those not satisfying this condition. Formally, we choose $\alpha_0<\pi$, and construct $\Th^{1}, \Th^{2}$:
\begin{equation}
\begin{split}
\Th^{1}&=\{K\in\Th:\ \alpha_K\leq\alpha_0\},\\
\Th^{2}&=\Th\setminus\Th^{1}.\label{th12}
\end{split}
\end{equation}
Since elements from $\Th^{1}$ satisfy the maximum angle condition, we can use Lagrange interpolation on them to obtain the element-wise $O(h)$ estimate implied by Lemma \ref{lem:th1est}. However, on $\Th^2$, Lemma \ref{lem:th1est} does not hold, therefore we modify the Lagrange interpolation conditions in this way: We seek $v_h^K\in P^1(K)$ such that
\begin{equation}
\begin{split}
\label{eq:th2Kdef}
\big(\grad\vhK(x_K),v_2\big)&=(\grad u(x_K),v_2),\\
\vhK(B_K)&=u(B_K),\\
\vhK(C_K)&=u(C_K).
\end{split}
\end{equation}
Here $v_1=\frac{C_K-B_K}{|C_K-B_K|}$ is the unit vector in the direction of the edge $(C_K,B_K)$, vector $v_2$ is the unit vector perpendicular to $v_1$ and $x_K$ is the foot of the altitude from vertex $A_K$ (cf. Figure \ref{fig:triang}). The first condition from (\ref{eq:th2Kdef}) says that $\vhK$ and $u$ have the same derivative in the direction $v_2$ at point $x_K$. The following estimates hold.

\begin{lemma}
\label{lem:th2est}
For $K\in\Th^2$ and $\vhK$ defined by (\ref{eq:th2Kdef}), we have
\begin{align}
\big|u-\vhK\big|_{H^1(K)}&\leq  \sqrt{13}h_K|u|_\di|K|^{1/2}, \label{lem:th2est_est1}\\
\big|u(A_K)-\vhK(A_K)\big|&\leq \big(|x_K-B_K||x_K-C_K|+|x_K-A_K|^2\big)|u|_\di.
\label{lem:th2est_est2}
\end{align}
\end{lemma}
\proof
\textbf{$\mathbf{H^1(K)}$-estimate}: Let $x\in K$. Since $v_1$ and $v_2$ are orthogonal vectors of unit length, we can write 
\begin{equation}
\begin{split}
\grad u(x)&=(\grad u(x),v_1)v_1+(\grad u(x),v_2)v_2,\\
\grad \vhK(x)&=\tfrac{u(C_K)-u(B_K)}{|C_K-B_K|}v_1+(\grad u(x_K),v_2)v_2.
\end{split}
\nonumber
\label{lem:th2est:1}
\end{equation}
The expression for $\grad \vhK$ is obtained from the fact that $(\grad \vhK(x),v_1)$, the coefficient at $v_1$, is the directional derivative of $\vhK$ in the direction $v_1$, which is $\tfrac{u(C_K)-u(B_K)}{|C_K-B_K|}$, due to (\ref{eq:th2Kdef}). Since $v_1\perp v_2$, 
\begin{equation}
\begin{split}
&\bigl|\grad(u-\vhK)(x)\bigr|^2\\ 
&=\!\Bigl(\!\bigl(\grad u(x),v_1\bigr)\!-\!\tfrac{u(C_K)-u(B_K)}{|C_K-B_K|}\Bigr)^2 \!\!\!+\!\Bigl(\!\bigl(\grad u(x),v_2\bigr)-\bigl(\grad u(x_K),v_2\bigr)\!\Bigr)^2\!\!=:(A)+(B).
\nonumber
\end{split}
\end{equation}
\emph{Estimate of (B):} By the multivariate Taylor expansion, we have
\begin{equation}
\frac{\pd u}{\pd x_j}(x)\!=\!\frac{\pd u}{\pd x_j}(x_K)+(x-x_K)\cdot\grad\frac{\pd u}{\pd x_j}(\xi) \Longrightarrow \Big|\frac{\pd u}{\pd x_j}(x) -\frac{\pd u}{\pd x_j}(x_K)\Big|\!\leq\! h_K\sqrt{2} |u|_\di,
\nonumber
\end{equation}
where $\xi$ lies on the line between $x$ and $x_K$. Therefore
\begin{equation}
\big|\grad u(x)-\grad u(x_K)\big|\leq\sqrt{2} h_K\sqrt{2} |u|_\di =2h_K|u|_\di.
\nonumber
\end{equation}
We can conclude that
\begin{equation}
(B)=(\grad u(x)-\grad u(x_K),v_2)^2\leq |\grad u(x)-\grad u(x_K)|^2|v_2|^2 \leq 4h_K^2|u|_\di^2.
\nonumber
\end{equation}

\noindent\emph{Estimate of (A):} We have
\begin{equation}
(A)=\Big(\underbrace{\big(\grad u(x)-\grad u(B_K),v_1\big)}_{(A_1)}+\underbrace{\big(\grad u(B_K),v_1\big)-\tfrac{u(C_K)-u(B_K)}{|C_K-B_K|}}_{(A_2)}\Big)^2.
\nonumber
\end{equation}
The term $(A_1)$ can be estimated similarly to $(B)$, hence $|(A_1)|\leq 2h_K|u|_\di$. As for $(A_2)$, by the multivariate Taylor expansion,
\begin{equation}
u(C_K)=u(B_K)+(C_K\!-\!B_K)\cdot\grad u(B_K)+\tfrac{1}{2}(C_K\!-\!B_K)^\T\grad^2 u(\xi)(C_K\!-\!B_K),
\nonumber
\end{equation}
where $\grad^2 u$ is the Hessian matrix of $u$. Therefore, by the definition of $v_1$,
\begin{equation}
(A_2)\!=\!\tfrac{\grad u(B_K)\cdot(C_K\!-\!B_K)}{|C_K-B_K|}-\tfrac{u(C_K)-u(B_K)}{|C_K-B_K|} =\tfrac{-1}{2|C_K-B_K|}(C_K\!-\!B_K)^\T\grad^2 u(\xi)(C_K-B_K).
\nonumber
\end{equation}
Thus $|(A_2)|\leq h_K|u|_\di$. Finally
\begin{equation}
\big|\grad(u-\vhK)(x)\big|^2=\big((A_1)+(A_2)\big)^2+(B)\leq 13h_K^2|u|_\di^2.
\nonumber
\end{equation}
Therefore, we can estimate the $H^1(K)$-error of $\vhK$:
\begin{equation}
\big|u-\vhK\big|_{H^1(K)}^2= \int_K \big|\grad(u-\vhK)(x)\big|^2\dx \leq 13h_K^2|u|_{2,\infty}^2|K|.
\nonumber
\end{equation}

\textbf{Vertex estimate}: On the edge $(B_K, C_K)$, the function $\vhK$ is simply the two-point linear Lagrange interpolation of $u$. For this we have the well known error expression $e(x)= \tfrac{1}{2}(x-x_0)(x-x_1)f^{\prime\prime}(\xi)$, where $f$ is the interpolated function and $\xi$ is between $x_0$ and $x_1$, cf. \cite{Davis}. In our case, $f$ is the restriction of $\vhK-u$ to the edge $(B_K,C_K)$ and $f^{\prime\prime}(\xi)=v_1^\T\grad^2(\vhK-u)(\xi) v_1=-v_1^\T\grad^2u(\xi) v_1$ is the second derivative $f$ in the direction $v_1$. Therefore
\begin{equation}
\begin{split}
\big|&\vhK(x_K)-u(x_K)\big|\\ &\leq\tfrac{1}{2}|x_K-B_K||x_K-C_K||v_1^\T\grad^2u(\xi) v_1|
\leq |x_K-B_K||x_K-C_K||u|_\di,
\label{lem:th2est:2}
\end{split}
\end{equation}
since $|v_1|=1$. By the multidimensional Taylor expansion, we have
\begin{equation}
\begin{split}
\big|u(A_K)&\!-\!\vhK(A_K)\big|\leq\big|u(x_K)\!-\!\vhK(x_K)\big| +\big|(\grad u(x_K)\!-\!\grad\vhK(x_K))\cdot(A_K\!-\!x_K)\big|\\
&+\big|\tfrac{1}{2}(A_K-x_K)^\T\grad^2 u(\xi)(A_K-x_K)\big|=(C_1)+(C_2)+(C_3).
\nonumber
\end{split}
\end{equation}
Term $(C_1)$ can be estimated by (\ref{lem:th2est:2}). For $(C_2)$, we have
\begin{equation}
\begin{split}
(C_2)&=\big|(\grad u(x_K)-\grad\vhK(x_K))\cdot\frac{(A_K-x_K)}{|A_K-x_K|}|A_K-x_K|\big|\\ &=\big|(\grad u(x_K)-\grad\vhK(x_K))\cdot v_2|A_K-x_K|\big|=0,
\nonumber
\end{split}
\end{equation}
due to the first condition in $(\ref{eq:th2Kdef})$. Finally, we have $(C_3)\leq |A_K-x_K|^2|u|_\di$. Collecting these estimates gives us (\ref{lem:th2est_est2}).
\qed

\begin{remark}
In (\ref{lem:th2est_est1}), instead of the expected $|u|_{H^2(K)}$, we have $|u|_\di|K|^{1/2}$. This term is an upper bound for $|u|_{H^2(K)}$ and mimics its behavior in that summing its squares over $K\in\Th$ gives $|u|_\di^2|\Omega|$, a constant. In this sense, (\ref{lem:th2est_est1}) is an $O(h_K)$ estimate.
\end{remark}

\begin{remark}
\label{rem:hk2}
A crude estimate of (\ref{lem:th2est_est2}) gives us $\big|u(A_K)-\vhK(A_K)\big|\leq 2h_K^2|u|_{2,\infty}$. In other words, we can fix the $O(h_K)$ interpolation property of the Lagrange interpolation in $H^1(K)$ if we make an $O(h_K^2)$ perturbation to the interpolated value of $u$ at $A_K$. At $B_K$ and $C_K$ we interpolate $u$ exactly. 
\end{remark}

If $\vhK$ is constructed by (\ref{eq:th2Kdef}), then $\vhK(A_K)\neq u(A_K)$, but by Lemma \ref{lem:th2est}, the difference of these values is $O(h_K^2)$. If $A_K$ also belongs to some element $K^\prime\in\Th^1$, then $v_h^{K^\prime}(A_K)=u(A_K)$, by the conditions $(\ref{eq:th1Kdef})$. Therefore, the global $V_h$-interpolation would not be continuous at $A_K$. We could fix the value $v_h^{K}(A_K)$ and modify $v_h^{K^\prime}$ so that $v_h^{K^\prime}(A_K)=\vhK(A_K)$. That would mean imposing more restrictive conditions on $K^\prime$ so that $v_h^{K^\prime}$ still satisfies estimate (\ref{lem:th1est:1}), even though we changed one of its vertex values by $O(h_K^2)$. For example, if the neighboring element $K^\prime$ was very small, e.g. $\diam K^\prime=O(h_K^2)$, then a change of vertex value of $O(h_K^2)$ would lead to an $O(1)$ change in $\grad v_h^{K^\prime}$. 

Another possibility is the following. We want to ``distribute" the $O(h_K^2)$ perturbation of the vertex value to a neighborhood of $A_K$, not only to the immediately neighboring element. If we can do this smoothly and locally, we can manage to preserve the $O(h)$ interpolation estimates on elements in that neighborhood. For this purpose we construct a so-called correction function for the modified Lagrange interpolation.

\begin{definition}
\label{def:w}
Let $\Omega_1=\cup_{K\in\Th^1}\overline K$ be the subset of $\Omega$ containing $\Th^1$. We call $w:\Omega\to\IR$ a \emph{correction function} corresponding to $u$ and $\Th$, if $w\in C(\overline\Omega)$, $w\in H^2(\Omega_1)$ and 
\begin{align}
(i)&\ w(x)=\vhK(x)-u(x) \text{ for all } K\in\Th^2 \text{ and all } x\in\{A_K, B_K, C_K\},\label{def:w1}\\
(ii)&\ |w|_{H^1(\Omega_1)}\leq C(u)h,\label{def:w2}\\
(iii)&\ |w|_{H^2(\Omega_1)}\leq C(u).\label{def:w3}
\end{align}
\end{definition}

\begin{remark}
In our context, the constants $C(u)$ in (\ref{def:w2}), (\ref{def:w3}) will take on the form $C(u)=C(|u|_2+|u|_{2,\infty})$, where $C$ is independent of $h,u$.
\end{remark}

\begin{remark}
Since $w\in C(\overline\Omega)$ and $(i)$ holds, for every $K\in\Th^2$ the maximum-angle vertex $A_K$ cannot be the vertex of any other $K^\prime\in\Th^2$, since $(i)$ would prescribe two different values of $w$ at this point. On the other hand, $K$ and $K^\prime$ can share the $B,C$ vertices, so we can have e.g. $B_K=C_{K^\prime}$. This is because $(i)$ prescribes $w(B_K)=w(C_{K^\prime})=0$ due to (\ref{eq:th2Kdef}). This is used in Lemma \ref{lem:w}.
\end{remark}

Now we shall proceed as follows. Instead of taking $v_h:=\Pih u$ in (\ref{Cea}), we shall take $v_h:=\Pih \tu,$ where $\tu=u+w$. Since $w\in C(\Omega)$, $\Pih\tu$ is well defined. Moreover, due to condition (\ref{def:w1}) for $K\in\Th^2$ and $x\in\{A_K, B_K, C_K\}$
\begin{equation}
\tu(x)=u(x)+\vhK(x)-u(x)=\vhK(x).
\nonumber
\end{equation}
Thus $(\Pih\tu)|_K=\vhK$ for all $K\in\Th^2$, hence $u-\Pih\tu$ satisfies the estimates of Lemma \ref{lem:th2est} on $\Th^2$. Finally, on $\Th^1$, where $\Pih u$ is a good approximation of $u$, $|w|_{H^1(\Omega_1)}=O(h)$, hence $\tu-u=O(h)$ in $H^1(\Omega)$, and something similar can be expected to hold for $\Pih\tu-\Pih u$.

\begin{lemma}
\label{lem:tu_interp}
Let $u\in\Wi$ and $w$ be as in Definition \ref{def:w}. Let $\tu=u+w$. Then 
\begin{equation}
|u-\Pih\tu|_1\leq C(u)h,
\nonumber
\end{equation}
where $C(u)$ is independent of $h$.
\end{lemma}
\proof
We have 
\begin{equation}
|u-\Pih\tu|_1^2=\sum_{K\in\Th^1}|u-\Pih\tu|_\HK^2 +\sum_{K\in\Th^2}|u-\Pih\tu|_\HK^2.
\label{lem:tu_interp1}
\nonumber
\end{equation}
Since $(\Pih\tu)|_K=\vhK$ for all $K\in\Th^2$, due to Lemma \ref{lem:th2est} we have
\begin{equation}
\sum_{K\in\Th^2}|u-\Pih\tu|_\HK^2\leq \sum_{K\in\Th^2}13h_K^2|u|_{2,\infty}^2|K|\leq 13h^2|\Omega||u|_{2,\infty}^2.
\label{lem:tu_interp2}
\nonumber
\end{equation}
On the other hand, due to (\ref{def:w2}), (\ref{def:w3}) and Lemma \ref{lem:th1est}
\begin{equation}
\begin{split}
&\sum_{K\in\Th^1}|u-\Pih\tu|_\HK^2\leq 2\sum_{K\in\Th^1}\Big(|u-\tu|_\HK^2+|\tu-\Pih\tu|_\HK^2\Big)\\ 
&\ \leq 2\!\!\sum_{K\in\Th^1}\Bigl(|w|_\HK^2+Ch_K^2|\tu|_\HdK^2\Bigr)\leq 2|w|_{H^1(\Omega_1)}^2+Ch^2(|u|_2^2+|w|_{H^2(\Omega_1)}^2)\\
&\ \leq (C(u)+|u|_2^2)h^2.
\nonumber
\label{lem:tu_interp3}
\end{split}
\end{equation}
We note that Lemma \ref{lem:th1est} can be applied since $\Th^1$ consists of elements satisfying the maximum angle condition, cf. (\ref{th12}). Combining the last three estimates gives us the desired result.
\qed

\subsection{Construction of the correction function $w$}
\label{subsec:constr_w}
\begin{figure}[t]
\begin{center}
\includegraphics[scale=0.6,clip]{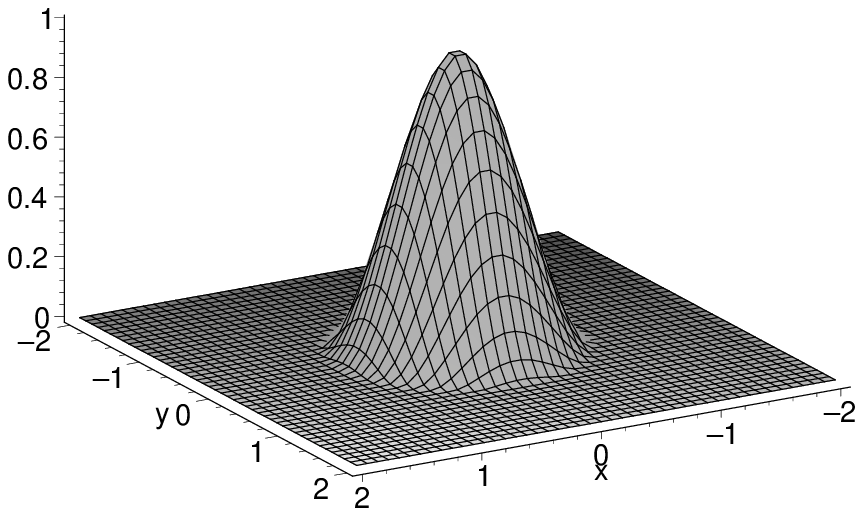}
\includegraphics[scale=0.6,clip]{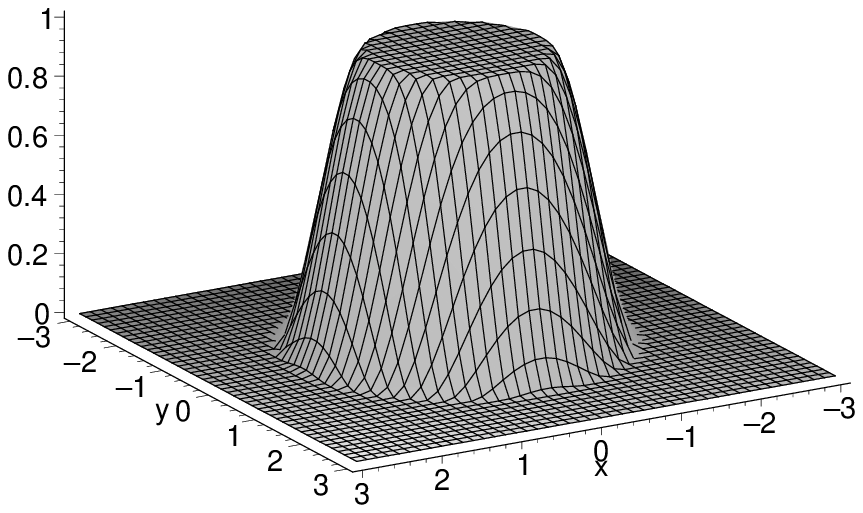}
\caption{Functions $\varphi_r$ and $\psi_r$ from (\ref{phi_r}) and (\ref{psi_r}) with $r=1$.}
\label{fig:phi_psi}
\end{center}
\end{figure}

We shall construct $w$ as a linear combination of disjoint local `bumps'. We define the cubic spline function $\tilde\vphi:[0,\infty)\to\IR$ as
\begin{equation}
\tilde\vphi=\begin{cases} 2(x-1)^3+3(x-1)^2,& \quad  x\in[0,1], \\
0,& \quad  x>1. \end{cases} 
\nonumber
\end{equation}
This function satisfies $\tilde\vphi\in C^1(0,\infty)\cap H^2(0,\infty)$ and $\tilde\vphi(0)=1, \tilde\vphi^\prime(0)=0, \tilde\vphi(1)=0, \tilde\vphi^\prime(1)=0$. Its derivatives are bounded by $|\tilde\vphi^\prime|\leq\tfrac{3}{2}, |\tilde\vphi^{\prime\prime}|\leq 6$. Using this function, we construct a local 2D bump with radius $r$:
\begin{equation}
\vphi_r(x)=\tilde\vphi(\tfrac{|x|}{r}),
\label{phi_r}
\end{equation}
cf. Figure \ref{fig:phi_psi}. We have $\vphi_r\in C^1(\IR^2)\cap H^2(\IR^2)$ and $\supp\vphi_r=\B_r(0)$, the disk with diameter $r$ centered at $0$.

\begin{lemma}
\label{lem:phi_r}
For all $x\in\IR^2$
\begin{equation}
\begin{split}
|\grad\vphi_r(x)|\leq\frac{3}{2r},\\
\|\grad^2\vphi_r(x)\|_F\leq\frac{36}{r^2},
\nonumber
\end{split}
\end{equation}
where $\|\cdot\|_F$ is the matrix Frobenius norm.
\end{lemma}
\proof
The presence of $1/r$ and $1/r^2$ can be seen from simple scaling arguments. The specific constants in the estimates can be obtained by a straightforward, yet somewhat tedious calculation.
\qed

\begin{figure}[t]
\begin{center}
\includegraphics[scale=0.58,clip]{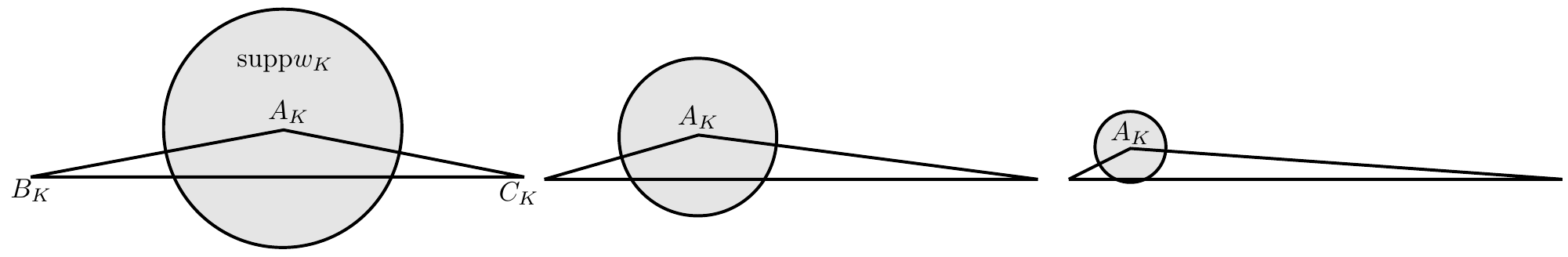}
\caption{Support of correction function $w_K$ for three different elements $K$ with the same maximum angle $\alpha_K$.}
\label{fig:phi_support}
\end{center}
\end{figure}

\medskip
First, we consider the case when $\Th^2$ consists of only one element $K$. Then $w$, corresponding to the single element $K$, can be constructed e.g. by
\begin{align}
w_K(x)&=\big(\vhK(A_K)-u(A_K)\big)\vphi_{r_K}(x-A_K),\label{eq:wKdef}\\
r_K&=\tfrac{1}{2}\min\{|A_K-B_K|,|A_K-C_K|\}.\label{eq:wKdef:rK}
\end{align}
The support of $w_K$ is shown in Figure \ref{fig:phi_support}. We note that the constant $\tfrac{1}{2}$ in (\ref{eq:wKdef:rK}) can be replaced by any $c\in(0,1]$, however we must keep in mind that this constant will appear in the denominator in estimates of Lemma \ref{lem:phi_r}. Finally, it is important to note that $B_K$ nor $C_K$ never lie in $\supp w_K$, since $r_K$ is half the distance from $A_K$ to the nearer of $B_K,C_K$.

\begin{lemma}
\label{lem:wK}
The function $w_K$ defined by (\ref{eq:wKdef}) is a correction function for $\Th^2$ consisting of one element $K$.
\end{lemma}
\proof
Since $\vphi_\rK(0)=1$, we have $w_K(A_K)=\vhK(A_K)-u(A_K)$. Furthermore, since $B_K,C_K\notin\supp\vphi_\rK(\cdot-A_K)$, we have $w_K(x)=0=\vhK(x)-u(x)$ for $x=B_K,C_K$, since in these two points $\vhK$ is an exact interpolation of $u$, due to (\ref{eq:th2Kdef}). Hence (\ref{def:w1}) is valid.

Since $w_K\in H^2(\Omega)$, the regularity condition is satisfied. It remains to prove (\ref{def:w2}), (\ref{def:w3}). Due to (\ref{eq:wKdef:rK}), either $|A_K-B_K|=2\rK$ and $|A_K-C_K|\leq h_K$ or vice versa. Moreover $|x_K-A_K|\leq 2r_K\leq h_K$. Hence, by (\ref{lem:th2est_est2})
\begin{equation}
\begin{split}
\big|\vhK(A_K)-u(A_K)\big|\leq\big(|x_K-B_K||x_K-C_K|+|x_K-A_K|^2\big) |u|_\di\\
\leq \big(|A_K-B_K||A_K-C_K|+|x_K-A_K|^2\big) |u|_\di\leq 4\rK h_K|u|_\di. 
\nonumber
\end{split}
\end{equation}
By Lemma \ref{lem:phi_r}
\begin{align}
|&w_K|_1^2=\!\int_{B_\rK\!(\!A_K\!)}|\grad w_K(x)|^2\dx \leq \int_{B_\rK\!(\!A_K\!)} \!\bigl|\vhK(A_K)-u(A_K)\bigr|^2\Bigl(\frac{3}{2\rK}\Bigr)^2\!\dx\nonumber\\
&\leq |B_\rK(A_K)| (4\rK h_K|u|_\di)^2 \Bigl(\frac{3}{2\rK}\Bigr)^2 \leq |B_\rK(A_K)| 36h_K^2|u|_\di^2.
\label{lem:wK:1}
\end{align}
Therefore $|w_K|_1\leq |\Omega|^{1/2}6h|u|_\di=C(u)h$. Similarly,
\begin{equation}
\begin{split}
|w_K|_2^2&=\int_{B_\rK\!(\!A_K\!)}\|\grad^2 w_K(x)\|_F^2\dx \leq \int_{B_\rK\!(\!A_K\!)}\bigl|\vhK(A_K)-u(A_K)\bigr|^2 \Bigl(\frac{36}{\rK^2}\Bigr)^2\!\dx\\
&\leq |B_\rK\!(\!A_K\!)| (4\rK h_K|u|_\di)^2 \Bigl(\frac{36}{\rK^2}\Bigr)^2 \leq \pi\rK^2 144^2\Bigl(\frac{h_K}{\rK}\Bigr)^2|u|_\di^2\\
&=\pi 144^2h_K^2|u|_\di^2.
\nonumber
\end{split}
\end{equation}
Therefore $|w_K|_2\leq \pi^{1/2}144h_K|u|_\di=C(u)h\leq C(u)$.
\qed

\medskip
In the proof of Lemma \ref{lem:wK}, we have used the very crude estimates $|B_\rK(A_K)|\linebreak\leq|\Omega|$ in the estimation of $|w_K|_1$ and $h_K\leq 1$  in the estimation of $|w_K|_2$. However, these quantities can be used in an additive way in the case that $\Th^2$ consists of multiple elements, where the correction function $w$ is simply the sum of all $w_K$. One must however take care that the resulting $w$ satisfies Definition \ref{def:w}. 

\begin{lemma}
\label{lem:w}
Let $\Th^2$ be such that 
\begin{equation}
\begin{split}
(i)\ &\sum_{K\in\Th^2}h_K^2\leq C,\\
(ii)\ &|A_K-A_{K^\prime}|\ge r_K+r_{K^\prime}\text{ for } K,K^\prime\in\Th^{2},\\
(iii)\ &|A_K-B_{K^\prime}|\ge r_K,|A_K-C_{K^\prime}|\ge r_K\text{ for } K,K^\prime\in\Th^{2},
\nonumber
\end{split}
\end{equation}
where $r_K, r_{K^\prime}$ are defined by (\ref{eq:wKdef:rK}). Then
\begin{equation}
w(x):=\sum_{K\in\Th^2} w_K(x)
\label{lem:w_def}
\end{equation}
is a correction function for $\Th^2$, where $w_K$ is defined by (\ref{eq:wKdef}).
\end{lemma}

\proof
Condition $(ii)$ merely states that $\supp w_K\cap\supp w_{K^\prime}=\emptyset$ and condition $(iii)$ states that no vertex $B_{K^\prime}$ or $C_{K^\prime}$ of another element $K^\prime$ is contained in $\supp w_K$. Together, $(ii)$ and $(iii)$ ensure that condition (\ref{def:w1}) is satisfied, since for any $A_K,B_K,C_K$ of any element $K\in\Th^2$, the sum (\ref{lem:w_def}) contains at most one nonzero term giving the correct value of $w$ in that particular point: $w(B_K)=w(C_K)=0$ and $w(A_K)=v_h^K(A_K)-u(A_K)$ for all $K\in\Th^2$.
 
As for  (\ref{def:w2}) and (\ref{def:w3}),
\begin{equation}
\begin{split}
|w|_1^2&=\!\sum_{K\in\Th^2}\int_{B_\rK\!(\!A_K\!)}\!|\grad w_K(x)|^2\dx 
\leq \!\sum_{K\in\Th^2}\!\!|B_\rK(A_K)| (4\rK h_K|u|_\di)^2 \Bigl(\frac{3}{2\rK}\Bigr)^2\\
&\leq |\Omega| 36h^2|u|_\di^2,
\nonumber
\end{split}
\end{equation}
similarly as in (\ref{lem:wK:1}). Therefore $|w|_1\leq |\Omega|^{1/2}6h|u|_\di=C(u)h$. Finally,
\begin{equation}
\begin{split}
|w|_2^2&=\!\!\sum_{K\in\Th^2}\int_{B_{\rK}\!(\!A_K\!)}\!\!\|\grad^2 w_K(x)\|_F^2\dx \leq\!\! \sum_{K\in\Th^2}\!\! |B_\rK(A_K)| (4\rK h_K|u|_\di)^2 \Bigl(\frac{36}{\rK^2}\Bigr)^{\!2}\\
&=\!\! \sum_{K\in\Th^2}\pi\rK^2 144^2\frac{h_K^2}{\rK^2}|u|_\di^2\leq\pi 144^2C|u|_\di^2,
\nonumber
\end{split}
\end{equation}
due to condition $(i)$. Therefore $|w|_2\leq \pi^{1/2}144C^{1/2}|u|_\di=C(u)$.

\qed

\begin{theorem}
\label{thm:suff_main}
Let $\Th=\Th^1\cup\Th^2$, where $\Th^1$ satisfies the maximum angle condition and $\Th^2$ satisfies the conditions of Lemma \ref{lem:w}. Then $|u-U|_1\leq C(u)h$ for all $u\in\Wi$.
\end{theorem}
\proof
The proof is a trivial consequence of Lemmas \ref{lem:w} and \ref{lem:tu_interp} and C\'ea's lemma.
\qed

\begin{figure}[t]
\begin{center}
\includegraphics[scale=0.57,clip]{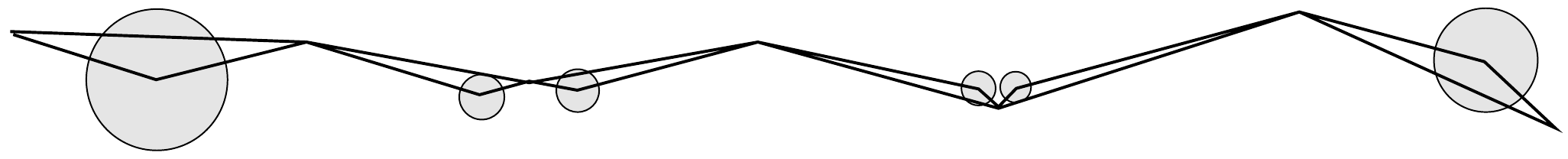}
\caption{Chain of degenerating elements with $\supp w_K$ depicted for each element.}
\label{fig:chains}
\end{center}
\end{figure}

\medskip
\noindent{\bf Example 1:} In Theorem \ref{thm:suff_main}, the maximum angle vertex $A_K$ of $K\in\Th^2$ must be sufficiently far from all other vertices of $\Th^2$, as in conditions $(ii),(iii)$. Since in $B_K,C_K$ the function $u$ is interpolated exactly, there is no condition at these vertices and elements from $\Th^2$ can be connected by these vertices. Still, one must be careful that $\supp w_K\cap\supp w_\tK=\emptyset$, and that $\supp w_K$ does not contain any other vertices in $\Th^2$. Once this is satisfied, elements from $\Th^2$ can be linked together to form arbitrarily long \emph{chains} such as in Figure \ref{fig:chains}. We note that the rest of the triangulation not indicated in Figure \ref{fig:chains}, i.e. $\Th^1$, can have an arbitrary structure, as long as it satisfies the maximum angle condition. Such triangulations then satisfy the assumptions of Theorem \ref{thm:suff_main}.

\begin{figure}[t]
\begin{center}
\includegraphics[scale=0.5,clip]{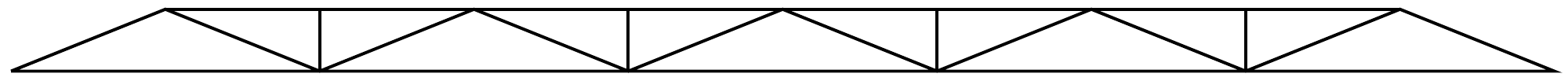}
\caption{Band with subdivided even elements satisfying Theorem \ref{thm:suff_main}.}
\label{fig:chainband}
\end{center}
\end{figure}

\medskip
\noindent{\bf Example 2:}
Theorem \ref{thm:suff_main} does not allow $\Th^2$ to contain the structures called bands in Section \ref{sec:necessary}, since they violate condition $(iii)$ of Lemma \ref{lem:w}, cf. Figure \ref{fig:Band}. However by subdividing the `even' elements of the band into two elements satisfying the maximum angle condition, we obtain the structure from Figure \ref{fig:chainband}, which satisfies the conditions of Theorem \ref{thm:suff_main}.

\begin{figure}[t]
\begin{center}
\includegraphics[scale=0.63,clip]{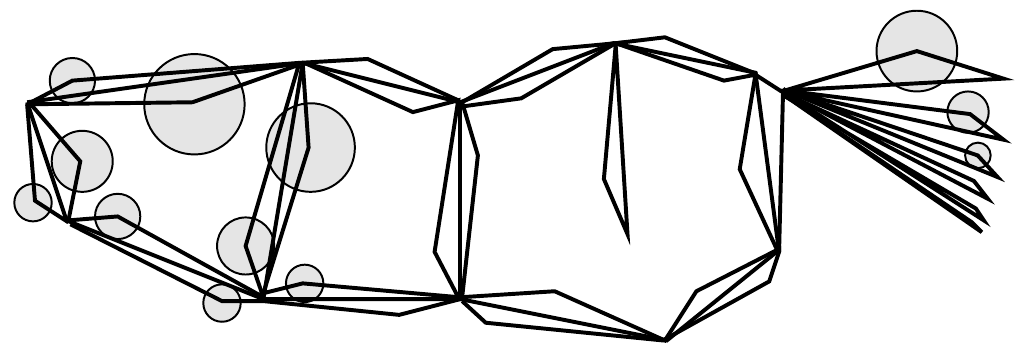}
\caption{More complex structures allowed in $\Th^2$ by Theorem \ref{thm:suff_main} with $\supp w_K$ depicted for some elements.}
\label{fig:chains2}
\end{center}
\end{figure}

\medskip
\noindent{\bf Example 3:}
Since linking of elements in $\Th^2$ is allowed as long as $A_K$ is not near another vertex from $\Th^2$ (i.e. $\supp w_K$ does not contain any other vertex from $\Th^2$), one can construct more complicated structures allowed by Theorem \ref{thm:suff_main} such as in Figure \ref{fig:chains2}. Again only the elements of $\Th^2$ are depicted, $\Th^1$ can fill in the remaining white spaces arbitrarily. Such structures can cover the whole domain $\Omega$, as long as condition $(i)$ of Lemma \ref{lem:w} is satisfied. As in the previous examples, by Theorem \ref{thm:suff_main}, the FEM has $O(h)$ convergence on such meshes. 

Obviously, one can fabricate very strange triangulations satisfying the assumptions of the theory presented. Such meshes perhaps do not have any value from the practical point of view, however when dealing with the mathematical question of necessary and sufficient conditions for convergence, one must consider general $\Th$, which can obviously be very `wild'.

\subsection{Clustering of elements}
\label{subsec:clusters}
\begin{figure}[t]
\begin{center}
\includegraphics[scale=0.7,clip]{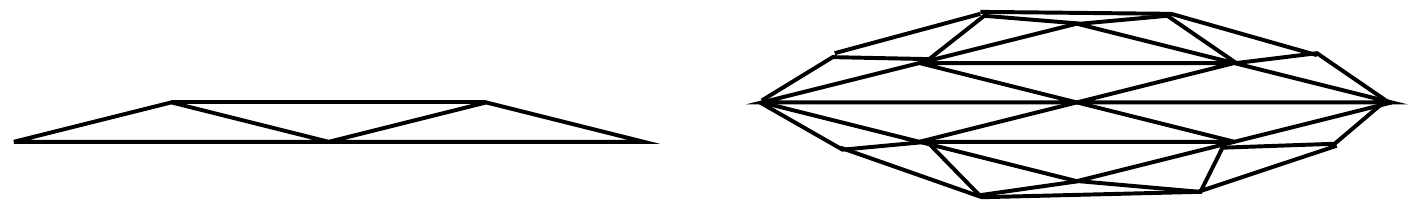}
\caption{Clusters consisting of degenerating elements.}
\label{fig:cluster}
\end{center}
\end{figure}

In Section \ref{subsec:modif_lagr}, we have modified the linear Lagrange interpolation on individual ``degenerating" elements $K\in\Th^2$. As we have seen, such elements can be connected to form ``chains", however their maximum angle vertices cannot be too close to other vertices from $\Th^2$, cf. Theorem \ref{thm:suff_main}. The question arises whether the presented construction can be extended also to situations such as Figure \ref{fig:cluster}, where the degenerating elements form nontrivial structures, which we will call \emph{clusters}, similar e.g. to the ``bands" of Section \ref{sec:necessary}. 

\begin{definition}
\label{def:cluster}
Let $\Th^1,\Th^2$ be as in (\ref{th12}). A \emph{cluster} of elements $\mC$ is a simply connected union of elements $K\in\Th$ such that there is at least one element $K\in\Th^2$ that lies in $\mC$.  
\end{definition}

\begin{remark}
We note that clusters, as from the previous definition, are not uniquely defined. Their choice for the purpose of proving $O(h)$ convergence is `user dependent' based on the specific geometry of $\Th$.
\end{remark}

In general, neither standard Lagrange interpolation (\ref{eq:th1Kdef}), nor the modified version (\ref{eq:th2Kdef}) can be used on elements from $\mC$. The first is not suitable for $K\in\Th^2$, the latter for $K,\tilde K\in\Th^2$ connected as in Figure \ref{fig:cluster}. The idea is to construct one linear function $\vhC$ globally on $\mC$ as the interpolation. Then we will use the idea of correction functions as in Section \ref{subsec:constr_w} to connect $\vhC$ continuously to the elements neighboring $\mC$. From $\vhC$ we require only that it is $O(h)$-close to $u|_\mC$ in the $H^1(\mC)$-seminorm with a constant independent of the geometry of $\mC$. The simplest possibility is to take $\vhC$ corresponding to the tangent plane of $u$ at some point $x_\mC\in\mC$.

Let $x_\mC\in\mC$ be an arbitrary but fixed point, we define $\vhC\in P^1(\mC)$ by
\begin{equation}
\vhC(x)=u(x_\mC)+\grad u(x_\mC)(x-x_\mC).
\label{eq:vhCdef}
\end{equation}
The following estimates are a straightforward consequence of Taylor's theorem.

\begin{lemma}
\label{lem:th2est_cluster}
Let $\vhC$ be defined by (\ref{eq:vhCdef}), then for $x\in\Omega$
\begin{align}
\big|u(x)-\vhC(x)\big|&\leq |x-x_\mC|^2|u|_\di,\nonumber\\ 
\big|\grad u(x)-\grad\vhC(x)\big|&\leq 2|x-x_\mC||u|_\di. \nonumber
\end{align}
\end{lemma}

\medskip
Similarly as in Definition \ref{def:w}, we define the correction function, which is used to continuously ``connect" the function $v_h^\mC$ to the rest of $\Th$ in the interpolation procedure.

\begin{definition}
\label{def:w_cluster}
Let $\{\mC_i\}_{i=1}^{N_\mC}$ be a set of clusters of elements from $\Th$. Let $\Omega_1=\Omega\setminus\cup_{i=1}^{N_\mC}\mC$. We call $w:\Omega\to\IR$ a \emph{correction function} corresponding to $u$ and $\Th$, if $w\in C(\overline\Omega)$, $w\in H^2(\Omega_1)$ and 
\begin{align}
(i)&\ w(x)=v_h^{\mC_i}(x)-u(x),\quad \forall i=1,\ldots N_\mC, \forall x\in\{A_K, B_K, C_K\},\forall K\in\mC_i,\label{def:w1_cluster}\\
(ii)&\ |w|_{H^1(\Omega_1)}\leq C(u)h,\label{def:w2_cluster}\\
(iii)&\ |w|_{H^2(\Omega_1)}\leq C(u).\label{def:w3_cluster}
\end{align}
\end{definition}

Similarly as in Section \ref{subsec:modif_lagr}, we shall use the correction function $w$ to construct a special interpolation of $u$ to use in C\'{e}a's lemma. We get the following theorem, the proof of which is essentially identical to that of Lemma \ref{lem:tu_interp}.

\begin{lemma}
\label{lem:tu_interp_cluster}
Let $u\in\Wi$ and $w$ be as in Definition \ref{def:w_cluster}. Let $\tu=u+w$. Then 
\begin{equation}
|u-\Pih\tu|_1\leq C(u)h,
\nonumber
\end{equation}
where $C(u)$ is independent of $h$.
\end{lemma}

As in Section \ref{subsec:constr_w}, we shall construct $w$ as a linear combination of disjoint local `bumps' around the individual clusters.  We define the spline function $\tilde\vphi:[0,\infty)\to\IR$ as
\begin{equation}
\tilde\psi=\begin{cases} 1,& \quad  x\in[0,1], \\ 2(x-2)^3+3(x-2)^2,& \quad  x\in[1,2], \\
0,& \quad  x>2. \end{cases} 
\nonumber
\end{equation}
This function satisfies $\tilde\vphi\in C^1(0,\infty)\cap H^2(0,\infty)$. Its derivatives are bounded by $|\tilde\psi^\prime|\leq\tfrac{3}{2}, |\tilde\psi^{\prime\prime}|\leq 6$. Using this function, we construct a local 2D `table mountain' bump with radius $2r$:
\begin{equation}
\psi_r(x)=\tilde\psi(\tfrac{|x|}{r}),
\label{psi_r}
\end{equation}
cf. Figure \ref{fig:phi_psi}. We have $\psi_r\in C^1(\IR^2)\cap H^2(\IR^2),\ \supp\psi_r=\B_{2r}(0)$ and $\psi_r\equiv 1$ on $B_r(0)$. Finally, we can estimate its derivatives similarly as in Lemma \ref{lem:phi_r}.

\begin{lemma}
\label{lem:psi_r}
For all $x\in\IR^2$
\begin{equation}
\begin{split}
|\grad\psi_r(x)|\leq\frac{3}{2r},\\
\|\grad^2\psi_r(x)\|_F\leq\frac{36}{r^2}.
\nonumber
\end{split}
\end{equation}
\end{lemma}

First, we consider the case when only one cluster $\mC$ is present in $\Th$. We define the corresponding correction function by
\begin{equation}
\begin{split}
w_\lmC(x)&=(v_h^\mC(x)-u(x))\psi_{r_\lmC}(x-x_\mC),\\
r_\lmC&=\diam\mC.
\end{split}
\label{eq:wCdef}
\end{equation}

\begin{lemma}
\label{lem:wC}
Let $\diam\mC\leq Ch^{1/2}$. The function $w_\lmC$ defined by (\ref{eq:wCdef}) is a correction function for $\Th$ containing one cluster $\mC$.
\end{lemma}
\proof
Since $\psi_r\equiv 1$ on $B_r(0)$, we have $\psi_r(x-x_\mC)=1$ for all $x\in\mC$, therefore $w_\mC(x)=v_h^\mC(x)-u(x)$ for all $x\in\mC$ Thus (\ref{def:w1_cluster}) is valid.

Trivially, $w_\lmC\in H^2(\Omega)$. It remains to prove (\ref{def:w2_cluster}), (\ref{def:w3_cluster}). By Lemmas \ref{lem:th2est_cluster} and \ref{lem:psi_r},
\begin{align}
&|w_\lmC|_1^2=\int_{\supp w_\lmC}|\grad w_\lmC|^2\dx  \label{lem:wC1}\\
&=\int_{\supp w_\lmC}\!\!\big|\grad (v_h^\mC(x)-u(x))\psi_{r_\lmC}\!(x-x_\mC) +(v_h^\mC(x)-u(x))\grad\psi_{r_\lmC}\!(x-x_\mC)\big|^2\!\dx\nonumber\\
&\leq \int_{\supp w_\lmC} \Big(4r_\mC|u|_{2,\infty}+4r_\mC^2|u|_{2,\infty}\frac{3}{2r_\mC}\Big)^2\dx \leq  400 r_\lmC^2|u|_{2,\infty}^2|\supp w_\lmC\!|\leq  Ch^2,
\nonumber
\end{align}
since $r_\mC^2=(\diam\mC)^2\leq Ch$ and $|\supp w_\lmC\!|=\pi(2r_\lmC)^2\leq Ch$. Hence $|w_\lmC|_1=O(h)$. Similarly (for brevity we omit the arguments of functions),
\begin{align}
&|w_\lmC|_2^2=\int_{\supp w_\lmC}\|\grad^2 w_\lmC\|^2_F\dx \nonumber\\
&=\int_{\supp w_\lmC}\big\|\grad^2 (v_h^\mC-u)\psi_{r_\lmC} +2\grad(v_h^\mC-u)\otimes\grad\psi_{r_\lmC}  +(v_h^\mC-u)\grad^2\psi_{r_\lmC}\big\|^2_F\dx\nonumber\\
&\leq \int_{\supp w_\lmC} \Big(2|u|_{2,\infty}+2r_\mC|u|_{2,\infty}\frac{3}{2r_\mC} +4r_\mC^2|u|_{2,\infty}\frac{36}{r_\mC^2}\Big)^2\dx\nonumber\\ 
&\leq C|\supp w_\lmC\!||u|_{2,\infty}^2\leq C(u),
\label{lem:wC2}
\end{align}
hence (\ref{def:w3_cluster}) is satisfied.
\qed

\subsubsection{Multiple clusters}
In the previous section, we have constructed a correction function for a single cluster $\mC$. If there are multiple clusters $\{\mC_i\}_{i=1}^{N_\mC}\subset\Th$, we can simply sum the individual correction functions for each $\mC_i$. In order to preserve the properties of the individual correction functions, we need to suppose their supports are disjoint. We get an analogy of Lemma \ref{lem:w}.

\begin{lemma}
\label{lem:w_cluster}
Let $\{\mC_i\}_{i=1}^{N_\mC}\subset\Th$ be such that 
\begin{equation}
\begin{split}
(i)\ &\sum_{i=1}^{N_\mC}r_{\mC_i}^2\leq Ch,\\
(ii)\ &\mathrm{dist}(\mC_i,\mC_j)\ge 2(r_{\mC_i}+r_{\mC_j})\ \text{ for } i,j=1,\dots,N_\mC,\ i\neq j,
\nonumber
\end{split}
\end{equation}
where $r_{\mC_i}=\diam \mC_i$. Then
\begin{equation}
w(x)=\sum_{i=1}^{N_\mC} w_{\lmC_i}(x)
\label{lem:w_def_cluster}
\end{equation}
is a correction function for $\Th$, where $w_{\mC_i}$ is defined by (\ref{eq:wCdef}).
\end{lemma}
\proof
As in Lemma \ref{lem:w}, condition $(ii)$ means that $\supp w_{\mC_i}\cap\supp w_{{\mC_j}}=\emptyset$ for $i\neq j$, i.e. property (\ref{def:w1_cluster}) is satisfied. As for  (\ref{def:w2}) and (\ref{def:w3}), due to (\ref{lem:wC1})
\begin{equation}
\begin{split}
|w|_1^2&=\!\sum_{i=1}^{N_\mC}\int_{\supp w_{\lmC_i}}\!\!\!|\grad w_{\lmC_i}(x)|^2\dx 
\leq \!\sum_{i=1}^{N_\mC}\!400 r_{\mC_i}^2|u|_{2,\infty}^2|\supp w_{\lmC_i}\!| \leq \!\sum_{i=1}^{N_\mC} C r_{\mC_i}^4|u|_{2,\infty}^2\\
&\leq C|u|_{2,\infty}^2\Big(\sum_{i=1}^{N_\mC} r_{\mC_i}^2\Big)^2\leq Ch^2|u|_{2,\infty}^2,
\nonumber
\end{split}
\end{equation}
due to assumption $(i)$. Therefore $|w|_1\leq C(u)h$. Similarly, we have due to (\ref{lem:wC2})
\begin{equation}
\begin{split}
|w|_2^2&=\sum_{i=1}^{N_\mC}\int_{\supp w_{\lmC_i}}\|\grad^2 w_{\lmC_i}(x)\|_F^2\dx \leq C|u|_{2,\infty}^2\sum_{i=1}^{N_\mC} |\supp w_{\lmC_i}\!| \leq C|u|_{2,\infty}^2|\Omega|,
\nonumber
\end{split}
\end{equation}
therefore $|w|_2\leq C(u)$.

\qed

We note that the constant 2 in condition $(ii)$ of Lemma \ref{lem:w_cluster} can be replaced by any $c>1$ by constructing the function $\tilde\psi$ such that $\psi\equiv 1$ on $[0,1]$ and $\supp \psi=[0,c]$. However we will get the factor $1/c$ in the corresponding version of Lemma \ref{lem:psi_r}.

\begin{theorem}
\label{thm:suff_main_cluster}
Let $\Th=\Th^1\cup\Th^2$, where $\Th^1$ satisfies the maximum angle condition and $\Th^2=\{\mC_i\}_{i=1}^{N_\mC}$ satisfies the conditions of Lemma \ref{lem:w_cluster}. Then $|u-U|_1\leq C(u)h$ for all $u\in\Wi$.
\end{theorem}
\proof
The proof is a trivial consequence of Lemmas \ref{lem:w_cluster} and \ref{lem:tu_interp_cluster} and C\'ea's lemma.
\qed

\medskip
\noindent{\bf Example 4:} Let $\{\mC_i\}_{i=1}^{N_\mC}\subset\Th$ be such that the cluster $\mC_i$ have diameter at most $Ch^{1/2}$ and are at least $Ch^{1/2}$ apart. Let $N_\mC\leq N$ for all $h\in(0,h_0)$, where $N$ is independent of $h$. Then the assumptions of Theorem \ref{thm:suff_main_cluster} are satisfied, hence the finite element method has $O(h)$ convergence.

\medskip
\noindent{\bf Example 5:} Let $\{\mC_i\}_{i=1}^{N_\mC}\subset\Th$ be such that the cluster $\mC_i$ have diameter at most $Ch$ and are at least $Ch$ apart. Let $N_\mC\leq Ch^{-1}$ for all $h\in(0,h_0)$. Then the finite element method has $O(h)$ convergence.

\medskip
\noindent{\bf Example 6:} Theorems \ref{thm:suff_main} and \ref{thm:suff_main_cluster} can be combined together by considering $\Th^2$ where clusters of elements of diameter at most $Ch^{1/2}$ exist together with the chains and other structures of Theorem \ref{thm:suff_main}. The corresponding correction function is simply the sum of all the particular correction functions. One only needs to ensure that the individual supports stay disjoint.

\medskip
\noindent We conclude with several remarks.

\begin{remark}
\label{rem:O_h_alpha}
Throughout Section \ref{sec:sufficient}, we have been interested in sufficient conditions for $O(h)$ convergence. If we want to generalise Theorems \ref{thm:suff_main} and \ref{thm:suff_main_cluster} to the case of $O(h^\alpha)$ convergence as in Section \ref{sec:necessary}, we can make the following changes:
\begin{itemize}
\item Everywhere, the maximum angle condition valid on $\Th^1$ can be replaced by the weaker condition: There exists $C_R>0$ such that $\Th^{1}=\{K\in\Th:\ R_K\leq C_Rh^\alpha\}$. By Lemma \ref{lem:circumradius}, this gives an $O(h^\alpha)$ estimate on $\Th^1$.
\item In condition (\ref{eq:wKdef:rK}), we can take $r_K=\tfrac{1}{2}h_K^{1-\alpha}\min\{|A_K-B_K|,|A_K-C_K|\}$, therefore, in condition $(ii)$ of Lemma \ref{lem:w}, the vertices $A_K,A_{K^\prime}$ can be much closer together for $\alpha<1$. On the other hand, condition $(i)$ becomes $\sum_{K\in\Th^2}h_K^{2\alpha}\leq C$, i.e. there must be fewer elements in $\Th^2$.
\item  On the other hand, we can keep the original choice of $r_K$ from (\ref{eq:wKdef:rK}). Thus we get an $O(h^\alpha)$ estimate on $\Th^1$ and $O(h)$ on $\Th^2$ with the original conditions $(i), (ii)$ of Lemma \ref{lem:w}.
\item In Lemma \ref{lem:wC}, we can have $\diam\mC\leq Ch^{\alpha/2}$, i.e the clusters can be much larger for $\alpha<1$. Consequently, in Lemma \ref{lem:w_cluster}, condition $(i)$ must be replaced by $\sum_{i=1}^{N_\mC}r_{\mC_i}^2\leq Ch^\alpha$. Then Theorem \ref{thm:suff_main_cluster} gives us $O(h^\alpha)$ convergence.
\end{itemize}
\end{remark}

\begin{remark}
The question arises whether the technique of Section \ref{sec:sufficient} is more general than the construction of \cite{Krizek}, where a `nice' triangulation satisfying the maximum angle condition is subdivided arbitrarily to still obtain $O(h)$ convergence on the subdivision. In other words, if $\Th$ satisfies the assumptions of Theorems \ref{thm:suff_main} and \ref{thm:suff_main_cluster}, does there exist a coarser triangulation $\tilde\Th$ with elements of diameter at most $Ch$ such that $\Th$ is a refinement of $\tilde\Th$? However in Theorems \ref{thm:suff_main} and \ref{thm:suff_main_cluster} nothing is assumed about $\Th^1$ except that it satisfies the maximum angle condition, it can have arbitrary structure. Then it is easy to construct $\Th$ containing e.g. one element in $\Th^2$, such that it is not a refinement of any triangulation. For example, if no two edges in $\Th$ sharing a vertex lie on a common line, then $\Th$ cannot be a refinement. 

Another deeper reason why $\Th$ in general is not a subdivision of a `nice' triangulation is the following lemma proved in \cite{Sbornik}:

\begin{lemma}
\label{lem:sbornik}
Let $\alpha\in(\tfrac{2}{3}\pi,\pi)$. Let $K$ be a triangle with all angles less than $\alpha$. Then there does not exist a finite conforming partition of $K$ into triangles which all contain an angle greater than or equal to $\alpha$.
\end{lemma}

In other words, if a triangle $K\in\tilde\Th$ is subdivided, the resulting subdivision must contain a triangle at least `as nice as' $K$. However Theorem \ref{thm:suff_main_cluster} allows for clusters of arbitrarily `bad' elements of size up to $Ch^{1/2}$, much larger than $Ch$ -- the size $K\in\tilde\Th$ would have to have in order to still have $O(h)$ convergence on $\tilde\Th$ and  the subdivision $\Th$.
\end{remark}

\begin{remark}
\label{rem:necandsuff}
The question of finding a necessary \emph{and} sufficient condition for $O(h^\alpha)$-convergence of the FEM remains open. However in some special cases the derived necessary (Section \ref{sec:necessary}) and sufficient (Section \ref{sec:sufficient}) conditions are not far apart. Take for example the situation in Counterexample 1 on page \pageref{nec:counter:1}. Let $\Th^2$ consist of the band $\B$ considered in this counterexample, i.e. with the shape parameter $\bar h=o(h^{4-8\alpha/5})$. Then for the length $L$ of $\B$ the necessary condition for $O(h^\alpha)$ convergence is $L\leq C_Lh^{2\alpha/5} =C_Lh^{0.4\alpha}$, while the sufficient condition is $L\leq C_Lh^{\alpha/2} =C_Lh^{0.5\alpha}$, since we can consider $\B$ as a cluster and apply Theorem \ref{thm:suff_main_cluster}. We note the difference in exponents is not large.

Similarly, if we consider $\Th^2$ consisting of multiple bands as in Counterexamples 6 and 7 (page \pageref{nec:counter:6}) forming a cluster of dimensions $L\times L$, then if $\bar h=o(h^{2-3\alpha/5})$ again the necessary condition for $O(h^\alpha)$ convergence is $L\leq C_Lh^{2\alpha/5}$, while the sufficient condition is $L\leq C_Lh^{\alpha/2}$.
\end{remark}

\section{Conclusion}
We have presented necessary and sufficient conditions on the triangulations $\Th$ for $O(h^\alpha)$ convergence of the piecewise linear conforming finite element method in $H^1(\Omega)$. The derived necessary condition limits the size of certain structures consisting of degenerating elements and is the first such condition for FEM convergence. The sufficient condition generalizes the maximum angle and circumradius conditions and shows that degenerating elements are allowed in $\Th$ as long as they form structures obeying certain rules. The analysis concerns the approximation properties of the discrete space $V_h$, since the finite element formulation itself is never used and only appears via C\' ea's lemma. Specifically:
\begin{itemize}
\item We have introduced the notion of bands, sets of neighboring elements connected by edges, such that their maximum angles form an alternating `zig-zag' pattern. Due to the geometry of the bands, the gradients of a continuous piecewise linear function $U$ on odd elements of the band determine the gradients on the intermediate even elements.

\item Provided $|u-U|_1=O(h^\alpha)$, this allowed us to produce a lower bound on the approximation error $|u-U|_{H^1(\B)}$ on the band $\B$, provided $\B$ is long enough ($L\ge C_L h^{2\alpha/5}$), cf. Theorem \ref{th:H1B_error}. 

\item The lower bound on $|u-U|_{H^1(\B)}$ can be made arbitrarily large e.g. by letting the maximal angles of $K\in\B$ go sufficiently fast to $\pi$ as $h\to 0$. Thus the basic assumption $|u-U|_1=O(h^\alpha)$ can be brought to a contradiction, giving necessary conditions on the geometry of $\B$ for the $O(h^\alpha)$ estimate to hold, cf. Corollary \ref{cor:H1B_error}.

\item Theorem \ref{th:mb:H1B_error} then gives the lower bound on $|u-U|_{H^1(\BB)}$ for a system of bands $\BB$, which gives the necessary condition for $O(h^\alpha)$ convergence of Corollary \ref{cor:mb:H1B_error}.

\item Based on the necessary conditions, various counterexamples to $O(h^\alpha)$ convergence were constructed, even examples of nonconvergence of the FEM. As a special case, we recovered the Babu\v{s}ka-Aziz counterexample, \cite{Babuska-Aziz}, specifically the optimal results of \cite{Oswald} thereon.

\item As for the sufficient conditions for $O(h^\alpha)$ convergence, we treated $\alpha=1$, the general case follows by simple arguments. We split $\Th$ into two parts -- $\Th^1$ satisfying and $\Th^2$ violating the maximum angle condition. To construct a suitable interpolant from $V_h$, on $\Th^1$ we use Lagrange interpolation and a modified Lagrange procedure on $\Th^2$. The modified variant has $O(h_K)$ approximation properties in $H^1(K)$ independent of the geometry of $K$, however $u$ is not interpolated exactly in the maximum angle vertex of $K$ but with an $O(h_K^2)$ perturbation, cf. Lemma \ref{lem:th2est}.

\item Due to the $O(h_K^2)$ perturbations at certain vertices of $\Th^2$, the interpolants on individual elements cannot be continuously connected. To construct a globally continuous interpolant, we introduced the concept of correction functions, which modify $u$ locally so that the Lagrange interpolation of the new function $\tilde u$ corresponds to the modified Lagrange interpolation of $u$ on $\Th^{2}$. We gave constructions of the correction functions in several cases, leading to Theorems \ref{thm:suff_main} and \ref{thm:suff_main_cluster}. Elements violating the maximum angle condition can form clusters of diameter up to $O(h^{1/2})$ or arbitrarily large `chains' as long as their maximal-angle vertices are not too close to other vertices in $\Th^2$. Examples of such triangulations were provided.
\end{itemize}

Although a condition for FEM convergence that would be both necessary and sufficient remains unknown, future work will be devoted to narrowing and perhaps closing the gap between the derived conditions at least in special cases, if not in general. The unifying idea behind the presented analysis is considering the size of sets of elements on which the interpolant can or cannot be one globally defined linear function. Such considerations led to conditions on the diameter -- $O(h^{0.5\alpha})$ being sufficient and $O(h^{0.4\alpha})$ being necessary, as in Remark \ref{rem:necandsuff}. Another subject for future work is to strengthen the results of Section \ref{sec:sufficient} to hold for $u\in H^2(\Omega)$, since currently $u\in W^{2,\infty}(\Omega)$ is needed.




\end{document}